\documentclass[12pt]{amsart}
\usepackage{amssymb,amsfonts,amsmath,amsopn,amstext,amscd,latexsym,xy,color,verbatim,MnSymbol}

\input xy
\xyoption{all}

\theoremstyle{remark}

\newtheorem{para}{\bf}[subsection]
\newtheorem{example}[para]{\bf Example}

\newtheorem{rem}[para]{\bf Remark}
\theoremstyle{definition}

\theoremstyle{plain}

\newtheorem{thm}[para]{Theorem}
\newtheorem{lemma}[para]{Lemma}
\newtheorem{cor}[para]{Corollary}
\newtheorem{prop}[para]{Proposition}

\newenvironment{numequation}{\addtocounter{para}{1}
\begin{equation}}{\end{equation}}

\DeclareMathOperator{\der}{der}

\DeclareMathOperator{\Gal}{Gal}

\DeclareMathOperator{\GL}{GL}
\DeclareMathOperator{\gr}{gr}

\DeclareMathOperator{\Hom}{Hom}

\DeclareMathOperator{\Ind}{Ind}

\DeclareMathOperator{\Lie}{Lie}

\DeclareMathOperator{\pr}{pr}
\DeclareMathOperator{\red}{red}

\DeclareMathOperator{\Res}{Res}

\DeclareMathOperator{\SL}{SL}

\newcommand{\eps}{{\epsilon}}
\newcommand{\vphi}{{\varphi}}

\newcommand{\vep}{{\varepsilon}}

\newcommand{\bbF}{{\mathbb F}}

\newcommand{\bbN}{{\mathbb N}}

\newcommand{\bbQ}{{\mathbb Q}}
\newcommand{\bbR}{{\mathbb R}}

\newcommand{\bbZ}{{\mathbb Z}}

\newcommand{\bB}{{\bf B}}

\newcommand{\bG}{{\bf G}}

\newcommand{\bM}{{\bf M}}

\newcommand{\bP}{{\bf P}}

\newcommand{\bS}{{\bf S}}
\newcommand{\bT}{{\bf T}}
\newcommand{\bU}{{\bf U}}

\newcommand{\bZ}{{\bf Z}}

\newcommand{\frb}{{\mathfrak b}}

\newcommand{\frg}{{\mathfrak g}}
\newcommand{\frh}{{\mathfrak h}}

\newcommand{\fro}{{\mathfrak o}}
\newcommand{\frp}{{\mathfrak p}}

\newcommand{\frs}{{\mathfrak s}}
\newcommand{\frt}{{\mathfrak t}}
\newcommand{\fru}{{\mathfrak u}}

\newcommand{\frx}{{\mathfrak x}}

\newcommand{\frz}{{\mathfrak z}}

\newcommand{\frG}{{\mathfrak G}}
\newcommand{\frH}{{\mathfrak H}}

\newcommand{\frM}{{\mathfrak M}}
\newcommand{\frN}{{\mathfrak N}}

\newcommand{\frP}{{\mathfrak P}}

\newcommand{\frS}{{\mathfrak S}}
\newcommand{\frT}{{\mathfrak T}}
\newcommand{\frU}{{\mathfrak U}}

\newcommand{\frX}{{\mathfrak X}}

\newcommand{\frZ}{{\mathfrak Z}}

\newcommand{\cA}{{\mathcal A}}

\newcommand{\cC}{{\mathcal C}}

\newcommand{\cF}{{\mathcal F}}
\newcommand{\cG}{{\mathcal G}}

\newcommand{\cR}{{\mathcal R}}

\newcommand{\cU}{{\mathcal U}}

\newcommand{\Qp}{{\mathbb Q_p}} 
\newcommand{\Cp}{\mathbb C_p} 
\newcommand{\Zp}{\mathbb Z_p} 

\newcommand{\qd}{\hfill $\square$}

\newcommand{\bksl}{\backslash}
\newcommand{\ra}{\rightarrow}
\newcommand{\lra}{\longrightarrow}
\newcommand{\hra}{\hookrightarrow}
\newcommand{\sub}{\subset}

\newcommand{\midc}{{\, | \,}}

\begin{document}

\title[Irreducibility of locally analytic principal series
representations]{On the irreducibility of locally analytic principal series
representations}
\author{Sascha Orlik}
\address{Institut f\"ur Mathematik, Universit\"at Paderborn,
Warburger Stra\ss{}e 100, D-33098 Paderborn,  Germany}
\email{orlik@math.uni-paderborn.de}
\author{Matthias Strauch}
\address{Indiana University, Department of Mathematics, Rawles Hall, Bloomington, IN 47405, USA}
\email{mstrauch@indiana.edu}
\thanks{M.S. is partially supported by NSF grant DMS-0902103.}
\maketitle

{\small {\bf Abstract.} Let $\bG$ be a $p$-adic connected
reductive group with Lie algebra $\frg$. For a parabolic subgroup
$\bP \sub \bG$ and a finite-dimensional locally analytic
representation $V$ of a Levi subgroup of $\bP$, we study the induced locally analytic
$\bG$-representation $W = \Ind_\bP^\bG(V)$. Our result is the
following criterion concerning the topological irreducibility of
$W$: if the Verma module $U(\frg) \otimes_{U(\frp)} V'$ associated
to the dual representation $V'$ is irreducible, then $W$ is
topologically irreducible as well.

\normalsize

\tableofcontents

\section{Introduction}

One of the principal methods for constructing representations of
reductive groups is to induce representations of parabolic
subgroups which come by inflation from representations of Levi
factors. This applies for example to the theory of algebraic
representations as well as to the theory of smooth representations
of $p$-adic reductive groups. In this paper we consider
parabolically induced representations in the theory of locally
analytic representations of $p$-adic reductive groups. A
systematic framework to study locally analytic representations of
$p$-adic groups was developed in the recent years mainly by P.
Schneider and J. Teitelbaum, cf. \cite{ST1}, \cite{ST2}. Algebraic
representations and smooth representations, as well as tensor
products of these, provide first examples of such locally analytic
representations, but there are much more. For instance, the
representations which are locally analytically induced from
representations of parabolic subgroups. Locally analytic principal
series representations for $\SL_2(L)$, $L$ a finite extension of
$\Qp$, were already defined and studied by Y. Morita in \cite{Mo}.
They were later reconsidered for the group $\GL_2$
in  \cite{ST1}  (for $L = \Qp$) and  in \cite{KS} (for arbitrary $L$).

\medskip

In this paper we prove a general criterion for the irreducibility
of parabolically induced locally analytic representations. In his
thesis H. Frommer \cite{Fr} studied locally analytic
representations of $\bG = \frG(\Qp)$, where $\frG$ is a split
reductive group over $\Qp$, which are induced from
finite-dimensional representations $V$ of a parabolic subgroup
$\bP \sub \bG$. $V$ is a $K$-vector space, where $K$ is a complete discretely valued field extension of $\Qp$.
His main theorem is a criterion for the
(topological) irreducibility of the induced representation
$\Ind^\bG_\bP(V)$ in terms of a canonically associated Verma
module. The crucial idea is to compare the structure of the dual
space $\Ind^\bG_\bP(V)'$  as a module over the distribution
algebra $D(G,K)$, where $G \sub \bG$ is a maximal compact
subgroup, with the structure of an associated Verma module over
the universal enveloping algebra $U(\frg)$
of the Lie algebra $\frg$ of $\bG$. In order to pass from
$D(G,K)$ to $U(\frg)$ one needs a technical result about their
relation and Frommer only showed this for $L = \Qp$. This is the
reason for the restriction to the base field $\Qp$ in \cite{Fr}. Later J.
Kohlhaase \cite{K1} proved this technical result for arbitrary
finite extensions $L$ over $\Qp$, which we use here to generalize
Frommer's theorem to the case of a not necessarily split reductive
group over an arbitrary extension of $\Qp$.

\medskip

In order to state our main result let $\frG$ be a connected
reductive group over $L$, and let $\frP \sub \frG$ be a parabolic
subgroup. We consider a locally analytic representation
$V$ of the group $\bP = \frP(L)$ which comes by inflation from a Levi
subgroup and put $\bG = \frG(L)$. Here $V$ is a vector space over a complete discretely
valued extension $K$ which contains $L$ and whose absolute value induces
the given absolute value on $L$. Then we have:

\medskip

{\bf Main result:} {\it Suppose $\dim_K(V) < \infty$. Then
the induced locally analytic
representation $\Ind^\bG_\bP(V)$ is topologically irreducible if
$U(\frg) \otimes_{U(\frp)} V'$ is irreducible as a module over the
universal enveloping algebra $U(\frg)$.}

\medskip

Our overall strategy of proof follows basically Frommer's
treatment. However, we found that an essential argument in
\cite{Fr}, stating that certain distribution algebras are integral
domains, is not obvious, as it is claimed there. We prove this
in sections \ref{Parahoric subgroups and their distribution algebras}
and \ref{appendix}.

\medskip

But beside the fact that it is desirable to have such a criterion
for the irreducibility in general, a motivation was provided by
the concrete example of certain $\Qp$-analytic principal series
representations of $\GL_2(\bbQ_{p^2})$ (regarded as a group over
$\Qp$). The interest in these representations
comes from conjectural relations to
two-dimensional crystalline representations of
$\Gal(\overline{\bbQ}_p|\bbQ_{p^2})$. Such
induced locally analytic representations play an
important role in the the $p$-adic Langlands
program, cf. \cite{BS}.

\medskip

In the last section we consider a particular example,
namely the case where $\frG$ comes by restriction of scalars
from a group which is split over $L$. This includes in particular the above mentioned case of the group
$\GL_2(\bbQ_{p^2})$ (regarded as a group over $\bbQ_p$).

\medskip

{\it Acknowledgements.} We would like to thank Tobias Schmidt and
Jan Kohlhaase for helpful discussions on distribution algebras. We
thank the SFB 478 ``Geometrische Strukturen in der Mathematik'' at
M\"unster for financial support of travel expenses. Finally, we
thank the referee for his careful reading, in particular for pointing
out a gap in a previous version. M.S. acknowledges support by the National Science Foundation under grant DMS-0902103.

\medskip

{\it Notation.} We let $p$ be a prime number and denote by $L$ a
finite extension of $\Qp$. The normalized $p$-adic (logarithmic)
valuation is denoted by $v_p$ (i.e. $v_p(p) = 1$). We let $K \subset \Cp$ be a
complete discretely valued field extension of $L$, whose absolute
value $|\cdot|_K$ is normalized so that $|p|_K = p^{-1}$. The rings of
integers are denoted by $\fro_L$ and $\fro_K$, respectively. For
notions and notation in the context of non-archimedean functional
analysis we refer to \cite{S}.

\bigskip

\section{Distribution algebras and locally analytic representations}

\subsection{Distribution algebras} In this section we recall some
definitions and results about algebras of distributions attached
to locally analytic groups, cf. \cite{ST1}, \cite{ST2}. We
consider a locally $L$-analytic group $H$ and denote by

$$C^{an}(H,K) = C^{an}_L(H,K)$$

\medskip

the locally convex $K$-vector space of locally $L$-analytic
functions on $H$ as defined in \cite{ST1}. The strong dual

$$D(H,K)=D_L(H,K):=(C^{an}_L(H,K))'_b$$

\medskip
is the topological algebra of $K$-valued distributions on $H$.
If furthermore $H$ is compact, then $D(H,K)$
has the structure of a Fr\'echet algebra. The multiplication
$\delta_1 \ast \delta_2$ of distributions $\delta_1, \delta_2 \in
D(H,K)$ is defined by

$$\delta_1 \ast \delta_2(f) =
(\delta_1 \otimes \, \delta_2)((h_1,h_2) \mapsto f(h_1h_2))
\,, $$

\medskip

where the distribution $\delta_1 \otimes \, \delta_2 \in D(H
\times H,K)$ has the property that for functions $f_1,f_2 \in
C^{an}(H,K)$, one has

$$(\delta_1 \otimes \, \delta_2)
((h_1,h_2) \mapsto f_1(h_1)f_2(h_2)) =
\delta_1(f_1) \delta_2(f_2) \,.$$

\medskip

The universal enveloping algebra $U(\frh)$ of the Lie algebra
$\frh = \Lie(H)$ of $H$ acts naturally on $C^{an}(H,K)$. On
elements $\frx \in \frh$, this action is given by

$$(\frx f)(h)
= \frac{d}{dt}(t \mapsto f(\exp(-t\frx)h))|_{t=0} \,.$$

\medskip

This gives rise to an embedding of
$U(\frh)_K := U(\frh) \otimes_L K$ into $D(H,K)$:

$$U(\frh)_K \hookrightarrow D(H,K) \,, \,\,
\frz \mapsto (f \mapsto (\dot{\frz}f)(1)) \,.$$

\medskip

Here $\frz \mapsto \dot{\frz}$ is the unique $K$-linear
anti-automorphism of $U(\frh)_K$ which induces multiplication by
$-1$ on $\frh$.

\bigskip

\subsection{Norms and completions of distribution algebras}
\label{Norms and completions}

\begin{para}\label{valuations} {\it $p$-valuations and global
charts.} Let $H$ be a compact locally $\Qp$-analytic group. Recall
that a map

$$\omega: H-\{1\} \lra (\frac{1}{p-1},\infty) \sub \bbR$$

\medskip

is called a {\it $p$-valuation} (cf. \cite{L}, III.2.1.2)
if the following conditions hold for all $g,h \in H$:\\

i) $\omega(gh^{-1}) \ge \min\{\omega(g), \omega(h)\}$,\\

ii) $\omega(g^{-1}h^{-1}gh) \ge  \omega(g) + \omega(h)$,\\

iii) $\omega(g^p) = \omega(g) + 1$.\\

\medskip

As usual one puts $\omega(1) = \infty$ and interprets the above
inequalities in the obvious sense, if a term $\omega(1)$ occurs.
Let $\omega$ be a $p$-valuation on $H$. The above conditions imply
that for any $\nu > 0$ the sets

$$H_\nu = \{h \in H \midc \omega(h) \ge \nu\} \mbox{ and }
H_{\nu^+} = \{g \in H \midc \omega(g) > \nu\}$$

\medskip

are normal subgroups of $H$. We put

$$\gr(H) = \bigoplus\nolimits_{\nu > 0} H_\nu / H_{\nu^+} \,.$$

\medskip

The commutator induces a Lie bracket on $\gr(H)$ which gives
$\gr(H)$ the structure of a Lie algebra over $\bbF_p$. The map
$\eps$ defined by

$$\eps: \gr(H) \ra \gr(H) \,, \,\, \eps(gH_{\nu^+}) = g^pH_{(\nu+1)^+}$$

\medskip

is an $\bbF_p$-linear map on $\gr(H)$, which gives $\gr(H)$ the
structure of a graded Lie algebra over $\bbF_p[\eps]$, cf.
\cite{L}, III.2.1.1. It is a free $\bbF_p[\eps]$-module, whose
rank is equal to the dimension of $H$ as a $\Qp$-analytic group,
cf. loc. cit., III.3.1.3/7/9.

\medskip

If $(h_1H_{\omega(h_1)^+}, \ldots, h_dH_{\omega(h_d)^+})$ is a
basis of $\gr(H)$ over $\bbF_p[\eps]$, then the elements $h_1, \ldots, h_d$
form a topological generating system of $H$, and the
following map

$$\bbZ_p^d \ra H \,, \,\, (a_1, \ldots, a_d) \mapsto
h_1^{a_1} \cdot \ldots \cdot h_d^{a_d}$$

\medskip

is well-defined and a homeomorphism. Moreover,

$$\omega(h_1^{a_1} \cdot \ldots \cdot h_d^{a_d})
= \min\{\omega(h_i) + v_p(a_i) \midc i = 1, \ldots, d \} \,.$$

\medskip

The sequence $(h_1,\ldots,h_d)$ is called a $p$-{\it basis} of $H.$
\end{para}

\medskip

\begin{para}\label{uniform groups} {\it Uniform pro-$p$ groups.}
We recall some definitions and results about pro-$p$ groups, cf.
\cite{DDMS}, ch. 3, 4. In this section $H$ will be a pro-$p$
group which is equipped with its topology as a pro-finite group.
Then $H$ is called {\it powerful} if $p$ is odd (resp. $p=2$) and
$H/\overline{H^p}$ (resp. $H/\overline{H^4}$ if $p = 2$) is
abelian. Here, $\overline{H^p}$ (resp. $\overline{H^4}$) is the
closure of the subgroup generated by the $p$-th (resp. fourth)
powers of its elements. If $H$ is topologically finitely generated
one can show that the subgroups $H^p$ (resp. $H^4$) are open and
hence automatically closed. The {\it lower $p$-series}
$(P_i(H))_{i \ge 1}$ of an arbitrary pro-$p$ group $H$ is defined
inductively by

$$P_1(H) = H \,, \,\,
P_{i+1}(H) = \overline{P_i(H)^p [P_i(H),H ]} \,.$$

\medskip

If $H$ is topologically finitely generated, then the groups
$P_i(H)$ are all open in $H$ and form a fundamental system of
neighborhoods of $1$, cf. Prop. 1.16 in loc.cit. A pro-$p$ group $H$ is
called {\it uniform} if it is topologically finitely generated,
powerful and its lower $p$-series satisfies

$$(H:P_2(H)) = (P_i(H):P_{i+1}(H))$$

\medskip

for all $i \ge 1$. If $H$ is a topologically finitely generated
powerful pro-$p$ group then $P_i(H)$ is a uniform pro-$p$ group
for all sufficiently large $i$, cf. loc.cit. 4.2. Moreover, any
compact $\Qp$-analytic group contains an open normal uniform
pro-$p$ subgroup, cf. loc.cit. 8.34.
\end{para}

\medskip

\begin{para}\label{canonical valuation} {\it The canonical
$p$-valuation on uniform groups.} Let $H$ be a uniform pro-$p$
group. It carries a distinguished $p$-valuation $\omega^{can}$
which is associated to the lower $p$-series and which we call the
{\it canonical $p$-valuation}. In order to define it, we let
$\vep_p = 2$ if $p = 2$ and $\vep_p = 1$ for odd $p$. For $h \neq
1$, we then put

$$\omega^{can}(h) = -1 + \vep_p + \max\{i \ge 1 \midc h \in P_i(H)\}
\,.$$

\medskip

To verify that this gives indeed a $p$-valuation one makes use of
the fact that \linebreak $[P_i(H),P_j(H)] \sub P_{i+j}(H)$ for all
$i,j \ge 1$, cf. Prop. 1.16 in  \cite{DDMS}. Also property (iii) follows from Prop. 2.7 in loc.cit.
For $p = 2$ one has to use the stronger statement that
$[P_i(H),P_j(H)] \sub P_{i+j+1}(H)$ for all $i,j \ge 1$, cf.
\cite{Sch}, proof of Prop. 2.1.
\end{para}

\medskip

We remark that a uniform pro-$p$ group $H$ has the property that the exponential map $\exp_{H}: \Lie(H) \dashrightarrow H$, which is, for a general $p$-adic Lie group, only defined on a sufficienly small lattice in $\Lie(H)$, is for uniform pro-$p$ groups defined on a unique $\Zp$-lattice $\Lambda \sub \Lie(H)$, which it maps bijectively onto $H$ (cf. \cite{DDMS}, sec. 4.5 and sec. 9.4).

\medskip

\begin{lemma}\label{L_uniformgroups} (i) Any compact locally $L$-analytic group $G$ has a normal open subgroup $H$ with the following properties:
\medskip

\begin{enumerate}
\item $H$ is uniform pro-$p$.

\medskip

\item Let $\Lambda \sub \Lie(H)$ be the $\Zp$-lattice which the exponential map $\exp_{H}: \Lie(H) \dashrightarrow H$ maps bijectively onto $H$; then $\Lambda$ is stable under multiplication by $\fro_L$, i.e., $\Lambda$ is an $\fro_L$-submodule of $\Lie(H)$.

\end{enumerate}

\medskip

(ii) Assume $H$ is uniform pro-$p$, and let $\Lambda$ be as in (2) above. Then, for any $\Zp$-basis $(\frx_1, \ldots, \frx_N)$ of $\Lambda$ the map

$$\Zp^N \ra H \,, \hskip6pt (a_i)_i \mapsto \exp_H(\frx_1)^{a_1} \cdot \ldots \cdot \exp_H(\frx_N)^{a_N}$$

\medskip

is a homeomorphism.

\medskip

\end{lemma}

{\it Proof.} (i) Let $\Lambda' \sub \Lie(G)$ be any $\fro_L$-lattice which is stable under the adjoint action of $G$ on $\Lie(G)$. For $t \gg 0$, we have that $[p^t\Lambda',p^t\Lambda'] \sub p^{\vep_p}p^t\Lambda'$. Hence $p^t\Lambda'$ is a uniform Lie algebra over $\Zp$ in the sense of \cite{DDMS}, sec. 9.4. For $t \gg 0$ the exponential function $\exp_G$ will be defined on $p^t\Lambda'$, and will map $p^t\Lambda'$ bijectively onto $H := \exp_G(p^t\Lambda')$. It follows that $H$ is then a uniform pro-$p$ group, cf. \cite{DDMS}, Thm. 9.10. $H$ is normal in $G$ because $\Lambda$ is invariant under the adjoint action of $G$.

\medskip

(ii) By \cite{DDMS}, Thm. 9.10, the $\Zp$-Lie algebra $L_H$ of $H$, as defined in \cite{DDMS}, sec. 4.5, is uniform. Denote by $\ast$ the multiplication on $L_H$ as defined in \cite{DDMS}, sec. 9.4, using the Baker-Campbell-Hausdorff series (which converges and is analytic on $L_H$, by \cite{DDMS}, Lemma 9.12).
Then $(L_H, \ast)$ is actually isomorphic to $H$, \cite{DDMS}, Thm. 9.10.
By \cite{DDMS}, Thm. 9.8, for any $\Zp$-basis $(\frx_i)$ of $L_H$ the corresponding elements of $H$ will be a topological generating set. \qed

\medskip

\begin{rem}\label{L_uniform} (i) Pro-$p$ groups $H$ which satisfy the the properties (1) and (2) of Lemma \ref{L_uniformgroups} will be considered repeatedly in this paper, and, for the purpose of this paper, we will call them {\it $L$-uniform}. $L$-uniform groups are exactly the groups that satisfy condition (L) in \cite{Sch}, cf. before Cor. 4.4 of \cite{Sch}, which is in turn equivalent to the conditions in Prop. 1.3.5 of \cite{K1}. We thought about using condition (L) as in \cite{Sch}, but found that the formulation above is slightly more natural.

\medskip

(ii) Let $H$ be a locally $L$-analytic uniform pro-$p$ group. Then the $\Zp$-Lie algebra $L_H$ of $H$ in the sense of \cite{DDMS}, sec. 4.5, is not necessarily an $\fro_L$-module. Even though one has an embedding $L_H \hra \Lie(H)$, and $\Lie(H)$ is an $L$-vector space, the $\Zp$-submodule $L_H$ need not be stable under multiplication by $\fro_L$. For a simple example, consider the abelian group $H = \Zp + p \fro_L$. The $\Zp$-Lie algebra $L_H$ of $H$ is $H$ itself, but $L_H$ is not an $\fro_L$-module if $L \neq \Qp$.
\end{rem}

\medskip

\begin{para}\label{norms} {\it Norms induced by $p$-valuations.}
In this section we let $H$ be a compact $L$-analytic group. The
distribution algebra $D(H,K)$ is then a {\it Fr\'echet-Stein
algebra} in the sense of \cite{ST2}, sec. 3. This means in
particular that there exists a family of norms $\|\; \|_r$,
$\frac{1}{p}<r<1$, on $D(H,K)$ such that, if $D_r(H,K)$ denotes
the completion of $D(H,K)$ with respect to $\|\; \|_r$, the
convolution product $*$ on $D(H,K)$ extends by continuity to a
product on $D_r(H,K)$. We recall briefly the construction of such
a family of norms. This is done in three steps, cf. \cite{ST2}
(proof of Thm. 5.1).

\bigskip

{\it Step 1.} We denote by $R^L_\Qp H$ the group $H$ when
considered as a locally $\Qp$-analytic group, and we put $d =
\dim(R^L_\Qp H)$. Let $H_0 \sub R^L_\Qp H$ be an open normal
subgroup which  is equipped with a $p$-valuation $\omega$.
For instance, by \cite{DDMS},
8.34, we may choose $H_0$ to be an open normal uniform pro-$p$
subgroup, and we may take for $\omega$ the canonical valuation
associated to its lower $p$-series (cf. \ref{canonical
valuation}). Then we consider the graded group $\gr(H_0)$ (which
depends on $\omega$). A choice of a basis of $\gr(H_0)$ gives rise
to a homeomorphism $\psi: \bbZ_p^d \ra H_0$, cf. \ref{valuations},
which in turn induces an isomorphism of locally convex $K$-vector
spaces

$$\psi^\ast: C^{an}(H_0,K) \stackrel{\simeq}{\lra} C^{an}(\bbZ_p^d,K)$$

\medskip

as well as an isomorphism of $K$-Banach spaces of continuous
functions

$$\psi^\ast: C(H_0,K) \stackrel{\simeq}{\lra} C(\bbZ_p^d,K) \,.$$

\medskip

Using Mahler expansions (\cite{L}, III.1.2.4) we can express
elements of $C(\bbZ_p^d,K)$ as series

$$f(x)=  \sum_{n \in \bbN_0^d} c_n {x \choose n} \,,$$

\medskip

where $c_n \in K$ and

$${x \choose n} := {x_1 \choose n_1} \cdot \ldots \cdot
{x_d \choose n_d}$$

\medskip

for the multi-indices $x = (x_1, \ldots, x_d)$ and $n = (n_1,
\ldots, n_d) \in \bbN^d_0$. Further, we have $|c_n| \rightarrow 0$
as $|n| = n_1 + \ldots + n_d \rightarrow \infty$. A continuous
function $f \in C(\bbZ_p^d,K)$ is locally analytic if and only if
$|c_n |r^{|n|} \rightarrow 0$ for some $r>1$. Consider the group
algebra $K[H_0]$ of $H_0$. By identifying elements of $H_0$ with
Dirac distributions, we get an embedding

$$K[H_0] \hookrightarrow D(H_0,K) \,.$$

\medskip

Put $b_i:=h_i-1 \in K[H_0]$ and set for $n \in \bbN^d_0,$

$${\bf b}^n = b_1^{n_1} \cdots b_d^{n_d} \,.$$

\medskip

Then we have ${\bf b}^n(f) = c_n$ for any continuous function $f
\in C(H_0,K)$, where the $c_n$ are the Mahler coefficients of
$\psi^*(f) \in C(\bbZ_p^d,K)$. It follows that every distribution
$\lambda \in D(H_0,K)$ has the shape

$$\lambda = \sum\nolimits_{n \in \bbN_0^d} d_n {\bf b}^n$$

\medskip

where $\{d_n r^{|n|} \midc n \in \bbN^d_0 \}$ is a bounded set for
all $0 < r < 1$. The norm $\| \cdot \|_r$ on $D(H_0,K)$ is then
defined by

$$\| \sum\nolimits_{n \in \bbN_0^d} d_n {\bf b}^n \|_r =
\sup \{|d_n|r^{\tau(n)}\midc n \in \bbN^d_0 \} \,.$$

\medskip

Here, $\tau(n)$ is given by $\tau(n) = \sum_i n_i \omega(h_i)$.
The Banach algebra $D_r(H_0,K)$ is defined to be the completion of
$D(H_0,K)$ with respect to $\| \cdot \|_r$. Thus we obtain

$$D_r(H_0,K) = \{ \sum\nolimits_{n \in \bbN_0^d}
d_n {\bf b}^n \in K[[b_1,\ldots,b_d]] \midc
\lim_{|n| \ra \infty} |d_n| r^{\tau(n)} = 0 \} \,.$$

\medskip

Furthermore, for $\frac{1}{p} < r < 1$, the norm $\| \cdot \|_r$
is multiplicative, cf. \cite{ST2} Thm. 4.5
and does not depend on
the chosen basis (loc. cit., before Thm. 4.11). We obtain a
projective system of noetherian Banach algebras such that

$$D(H_0,K)= \varprojlim\nolimits_r D_r(H_0,K) \,.$$

\medskip

Moreover, the transition maps

$$D_{r'}(H_0,K) \rightarrow D_r(H_0,K)$$

\medskip

are flat for $\frac{1}{p} < r \le r' < 1$.

\medskip

{\it Step 2.} We extend the norm $\| \cdot \|_r$ on $D(H_0,K)$ to
a norm on $D(R^L_\Qp H,K)$ as follows. Let $\eta_1, \ldots,
\eta_s$ be a system of coset representatives of $H/H_0$. Then the
Dirac distributions $\delta_{\eta_1}, \ldots, \delta_{\eta_s}$
form a basis of $D(R^L_\Qp H,K)$ over $D(H_0,K)$. Writing $\mu \in
D(R^L_\Qp H,K)$ as $\mu = \lambda_1 \delta_{\eta_1} + \ldots +
\lambda_s \delta_{\eta_s}$, we put

$$q_r(\mu) = \max\{\| \lambda_i \|_r \midc 1 \le i \le s \} \,.$$

\medskip

Again, we obtain a system of sub-multiplicative norms $q_r$,
$\frac{1}{p}<r<1$, and a projective system $D_r(R^L_\Qp H,K)$ of
$K$-Banach algebras fulfilling the conditions above, and we have
$D(R^L_\Qp H,K) = \varprojlim D_r(R^L_\Qp H,K)$. The norms $q_r$
do not depend on the chosen representatives. It is worth to remark
at this point that the norms $q_r$ are in general not
multiplicative, because $H$ may contain non-trivial elements of
finite order which causes the group ring, as well as the rings of
distributions, to have non-trivial zero divisors\footnote{In
\cite{Fr} it is mistakenly stated that these norms were
multiplicative so that the corresponding completions would be
integral domains. As a consequence, some completed distribution
algebras were claimed to
be obviously integral domains. 
}.

\medskip

{\it Step 3.} Locally $L$-analytic functions on $H$ are obviously
locally $\Qp$-analytic functions on $R^L_\Qp H$, and hence there
is a canonical map

$$C^{an}_L(H,K) \hookrightarrow C^{an}_\Qp (R^L_\Qp H,K)$$

\medskip

which is a closed embedding (cf. the proof of Thm. 5.1 in
\cite{ST2}). This map induces by duality a continuous surjection

$$D_\Qp (R^L_\Qp H,K) \ra D_L(H,K)$$

\medskip

on the distribution algebras. We denote the induced residue norm
on $D_L(H,K)$ respectively on the completion $D_r(H,K)$ by
$\bar{q}_r$. Once again, we have $D_L(H,K)= \varprojlim D_r(H,K)$.

\medskip

Although we will mostly consider the case where the parameter $r$
lies strictly between $\frac{1}{p}$ and $1$, we point out that the
norms $\bar{q}_r$ are defined, as norms of $K$-vector spaces,
for any $r \in (0,1)$. This is used later in the proof of
Prop. \ref{from_Mr to mr}, cf. Lemma \ref{density}.
\end{para}

\medskip

\begin{rem}\label{equivalent_definition} We keep the notation
from the preceding paragraph. In the recipe for the norm
$\bar{q}_r$ on $D_L(H,K)$ given above, we could have carried out
steps two and three in reverse order with the same resulting norm. To
be precise, let $\overline{\| \cdot \|}_r$ be the quotient norm on
$D_L(H_0,K)$ induced by the surjection $D_{\Qp}(R^L_\Qp H_0,K)
\twoheadrightarrow D_L(H_0,K)$. We have again

$$D_L(H,K) = \bigoplus_{i = 1}^s D_L(H_0,K) \ast \delta_{\eta_i}$$

\medskip

and can hence consider the maximum norm $\hat{q}_r$ on $D_L(H,K)$
induced from this decomposition and the norm $\overline{\| \cdot
\|}_r$ on $D_L(H_0,K)$. Then it is not difficult to check that
$\hat{q}_r$ coincides with the norm $\bar{q}_r$ defined above, cf.
\cite{Sch}, Lemma 4.4.
\end{rem}

\bigskip

\subsection{The closure of the enveloping algebra}

\begin{para}\label{notation closure}
{\it Definition of $U_r(\frh,H_0)$.} Let $H$ be a compact locally
$L$-analytic group and let $H_0 \sub H$ be an open normal subgroup,
which is equipped with a $p$-valuation $\omega$. Associated to
$\omega$ there is a norm $\bar{q}_r$ on $D(H,K)$ as defined in
\ref{norms}. We write

$$U_r(\frh, H_0) \subset D_r(H_0,K)$$

\medskip

for the topological closure of $U(\frh) \otimes_L K$ in
$D_r(H_0,K)$. The following theorem by Kohlhaase generalizes a
result of Frommer who considered the case of $\Qp$-analytic
groups. It is this basic technical result which we use to
generalize Frommer's irreducibility criterion.
\end{para}

\medskip

\begin{thm}\label{Kohlhaase} (\cite{K1}, 1.4.2.)
Let $H$ be a compact locally $L$-analytic group of dimension $d$
with Lie algebra $\frh = \Lie(H)$ and $r\in p^\bbQ$ with
$\frac{1}{p}<r<1$.

\medskip

(i) Let $H_0 \sub H$ be a normal open subgroup which is $L$-uniform (cf. Remark \ref{L_uniform}).
Then, if we equip $R^L_\Qp H_0$
with its canonical valuation (\ref{canonical valuation}),
$D_r(H_0,K)$ is a free, finitely generated module over the
noetherian subalgebra $U_r(\frh,H_0)$.

\medskip

(ii) Let the normal $L$-uniform subgroup $H_0 = \exp_H(\Lambda) \sub H$ be as above, and let
$\frX = (\frX_1, \ldots, \frX_d)$ be an $\fro_L$-basis of $\Lambda$. Then there is a norm $\nu_r$ on $U_r(\frh,H_0)$ which is
equivalent to $\bar{q}_r$, such that $U_r(\frh,H_0)$ consists of
exactly those series

$$\sum\nolimits_{n \in \bbN_0^d} d_n \frX^n$$

\medskip

for which

$$\lim_{|n| \ra \infty} |d_n| \nu_r(\frX^n) = 0 \,,$$

\medskip

and $\nu_r(\sum_n d_n \frX^n) = \sup_n |d_n| \nu_r(\frX^n)$.
\end{thm}

\medskip

\begin{rem} In \cite{K1}, 1.4.2, the conditions imposed on $H_0$ are those of \cite{K1}, 1.3.5. However, a group $H_0$ which is $L$-uniform automatically satisfies the conditions of \cite{K1}, 1.3.5.\footnote{On the other hand, one can show that a group satisfying the conditions of \cite{K1}, 1.3.5, is actually $L$-uniform. We will not need this fact in this paper.}
\end{rem}

\medskip

\begin{cor}\label{nontrivial intersection}
Let $H$ be a compact locally $L$-analytic group, and let $H_0 \sub
H$ be an open normal subgroup as in Thm. \ref{Kohlhaase} (i). Suppose
that $D_r(H,K)$ is an integral domain. Then any non-zero left
ideal $I$ of $D_r(H,K)$ has non-zero intersection with
$U_r(\frh,H_0)$.
\end{cor}

{\it Proof.} Compare \cite{Fr}, Corollary 5. Since $H$ is compact,
$D_r(H,K)$ is a finite $D_r(H_0,K)$-module, and so $D_r(H,K)$ is
finite over $U_r(\frh,H_0)$. Therefore, there exists for any $F
\in D_r(H,K),$ a polynomial $P = \sum_{i=0}^N a_iX^i$, $a_i \in
U_r(\frh,H_0)$ with $P(F) = 0$. Because we assumed $D_r(H,K)$ to
be an integral domain, we may find a polynomial with $a_0 \neq 0$.
So, if $F \in I$ then
$$0\neq a_0 = - \sum_{i=1}^N a_i F^i \in I \cap U_r(\frh,H_0) \,.$$
\qd

\medskip

\subsection{Locally analytic representations}\label{Locally analytic
representations}

We conclude this section by recalling some facts of locally
analytic representations. Let $H$ be a locally $L$-analytic group,
$V$ a Hausdorff locally convex $K$-vector space, and $\rho: H \ra \GL_K(V)$
a homomorphism. Then  $\rho$ (or the pair $(V,\rho)$) is called a
{\it locally analytic representation} of $H$
if the topological $K$-vector space $V$ is barrelled,
each $h \in H$ acts $K$-linearly and continuously on $V$, and the
orbit maps $\rho_v: H \rightarrow V, \; h \mapsto \rho(h)(v)$, are
locally analytic maps for all $v \in V$, cf. \cite{ST1}, sec. 3.
If $V$ is of compact
type, i.e., a compact inductive limit of Banach spaces, the strong
dual $V'_b$ is a nuclear Fr\'echet space and a separately continuous left
$D(H,K)$-module. The module structure is given as follows:

$$D(H,K) \otimes_K V'_b \ra V'_b \,, \,\,
\delta \otimes \vphi \mapsto (v \mapsto \delta(g \mapsto
\vphi(\rho(g^{-1})v))) \,.$$

\medskip

This functor gives an equivalence of categories

$$\begin{array}{ccc}

\left\{\begin{array}{c}
\mbox{locally analytic $H$-represen-} \\
\mbox{tations on $K$-vector spaces } \\
\mbox{of compact type with } \\
\mbox{continuous linear $H$-maps }
\end{array} \right\}

& \longrightarrow &

\left\{\begin{array}{c}
\mbox{separately continuous $D(H,K)$-}  \\
\mbox{modules on nuclear Fr\'echet spa-} \\
\mbox{ces with continuous $D(H,K)$- } \\
\mbox{module maps}
\end{array} \right\}

\end{array}
$$

\medskip

In particular, $V$ is topologically irreducible if $V'_b$
is a simple $D(H,K)$-module.

\bigskip
For any closed subgroup $H'$ of $H$ and any locally analytic
representation $V$ of $H',$ we denote by $\Ind^H_{H'}(V)$ the
induced locally analytic representation. We recall the definition:

\begin{eqnarray*}
\Ind^H_{H'}(V)& := &\Big\{f \in C^{an}(H,V) \midc f(h' \cdot
h) = h' \cdot f(h) \;\forall h' \in H', \forall h \in H
\Big\}.
\end{eqnarray*}

\medskip

The group $H$ acts on this vector space by $(h \cdot f)(x) = f(xh)$.
We have a Frobenius reciprocity in the category of locally analytic
representations (see \cite{Fe} Theorem 4.2.6):

$$\Hom^{\rm cont}_H(W,\Ind^H_{H'}V) \cong \Hom^{\rm cont}_{H'}(\Res^H_{H'}W,V) \, .$$

\medskip

Here, $\Res^{H}_{H'}(W)$ denotes as usual the restriction, viewing
$W$ via the embedding $H'\hookrightarrow H$ as a
$H'$-representation. Consider the canonical map of $D(H,K)$-modules

\begin{numequation}\label{dual_iso}
D(H,K) \otimes_{D(H',K)} V' \ra (\Ind^H_{H'}V)'_b \,, \,\,
\delta \otimes \vphi \mapsto \delta \cdot \vphi \,,
\end{numequation}

\medskip

with $(\delta \cdot \vphi)(f) = \delta(g \mapsto \vphi(f(g^{-1})))$. Suppose $H = C \cdot H'$ with a compact subgroup $C \sub H$, such that $C \cap H'$ is topologically finitely generated, and suppose, moreover, that $V$ is finite-dimensional. Then the same arguments as in \cite{ST3}, before Lemma 6.1, show that (\ref{dual_iso}) is an isomorphism of topological vector spaces, if we give the left side the quotient topology of the projective tensor product topology on $D(H,K) \otimes_K V'$.\footnote{We remark that the projective and inductive tensor product topologies coincide for tensor products of Fr\'echet spaces, cf. \cite{S}, Prop. 17.6.} We will only consider induced representations where these conditions are met. (For more general situations cf. \cite{K2}, Prop. 5.3 and Remark 5.4.)

\bigskip

\section{Representations induced from a parabolic subgroup}

\subsection{The setting and statement of the main result}

\begin{para}\label{the groups} We consider a connected reductive
algebraic group $\frG$ over the finite extension $L$ of $\Qp$. Let
$\frS \sub \frG$ be a maximal split torus over $L$. Fix a minimal
parabolic subgroup $\frP_0$ of $\frG$ which contains $\frS$ and
denote by $\frU_0$ its unipotent radical. The choice of $\frP_0$
determines subsets $\Phi^+ \supset \Delta$ of positive and simple
roots, respectively, in the root system $\Phi$ of $\frG$ with
respect to $\frS$. Let $\frP$ be a parabolic subgroup containing
$\frP_0$ with unipotent radical $\frU.$ Let $\frM$ be the Levi
subgroup of $\frP$ containing $\frS$. Let $W$ be the Weyl group of
$\frG$ with respect to $\frS$, and let $W_{\bP} \subset W$ be the Weyl
group of the Levi subgroup $\frM.$ Finally, we denote by
$\Phi^{\red}$ the set of reduced roots of $\Phi $.

\bigskip

We denote the corresponding groups of $L$-valued points by bold
letters:

$$\bG = \frG(L) \,, \,\, \bP = \frP(L) \,, \,\,
\bS = \frS(L) \,, \,\, \mbox{ etc. }$$

\medskip

which we consider as locally $L$-analytic groups. The Lie algebras
will be denoted by gothic letters, i.e.,

$$\frg = \Lie(\bG) \,, \,\, \frp = \Lie(\bP) \,, \,\,
\frs = \Lie(\bS) \,, \,\, \mbox{ etc. }$$
\end{para}

\begin{para}\label{induced representation} We let

$$\rho: \bP  \lra \bM \lra \GL(V)$$

\medskip

be a representation of $\bP$ which comes by inflation from a
locally $L$-analytic representation of $\bM$ on a {\it
finite-dimensional} $K$-vector space $V$. We are interested in the
locally analytic induced representation

$$\Ind^{\bG}_{\bP}(\rho)$$

\medskip

and our main result gives a criterion for the topological
irreducibility of this locally $L$-analytic representation of
$\bG$ in terms of the generalized Verma module

$$m(\rho) := U(\frg) \otimes_{U(\frp)} \rho' \,.$$

\medskip

Here $\rho'$ is the derived representation of $\frp$ on the dual
space $V' = \Hom_K(V,K)$. Our main result then is
\end{para}

\medskip

\begin{thm}\label{main result}
If $m(\rho)$ is a simple $U(\frg)$-module, then
$\Ind^{\bG}_{\bP}(\rho)$ is a topologically irreducible
representation.
\end{thm}

\bigskip

\subsection{The structure of the proof}\label{Overview of the proofI}

\begin{para} {\it Reduction to groups of the type (torus) $\times$ (simply connected).} The first step is to reduce to the case where $\frG$ is a product of a torus and a
semi-simple simply connected algebraic group.\footnote{This is needed to ensure that the Iwahori subgroup $I$, to be defined below, possesses Iwahori decompositions w.r.t. all parabolic subgroups $w\bP w^{-1}$, cf. \ref{defining parahoric subgroups}.} To this end let
$\tilde{\frG}_{\der}$ be the simply connected cover of the derived group
$\frG_{\der}$ of $\frG$. Let $\frZ$ be the connected component of the center of $\frG$. Put $\tilde{\frG} = \frZ \times \tilde{\frG}_{\der}$. There is a canonical morphism $\iota: \tilde{\frG} \ra \frG$ which induces an isomorphism between the corresponding Lie algebras. Put $\tilde{\frP} = \iota^{-1}(\frP)$, which is a parabolic subgroup of $\tilde{\frG}$. The morphism $\iota$ indices an isomorphism of flag varieties $\tilde{\frG}/\tilde{\frP} \stackrel{\simeq}{\lra} \frG/\frP$, from which we deduce an isomorphism of representations of $\tilde{\bG}$

$$\Ind^{\bG}_{\bP}(\rho) \stackrel{\simeq}{\lra} \Ind^{\tilde{\bG}}_{\tilde{\bP}}(\tilde{\rho})\,.$$

\medskip

Here $\tilde{\rho}$ is the representation of $\tilde{\bP}$ given
by the composition of the induced map $\tilde{\bP} \ra \bP$ and $\rho$. The action of $\tilde{\bG}$ on the right comes
from the homomorphism $\tilde{\bG} \ra \bG$. It follows that
$\Ind^{\bG}_{\bP}(\rho)$ is a (topologically) irreducible
$\bG$-representation if
$\Ind^{\tilde{\bG}}_{\tilde{\bP}}(\tilde{\rho})$ is a
(topologically) irreducible $\tilde{\bG}$-representation. The same
argument applies to the Verma modules. Therefore, we will assume
from now on that $\frG$ is the product of a central torus and a semi-simple and simply connected group.
\end{para}

\medskip

\begin{para}\label{Passage to representations of compact groups}
{\it Passage to representations of compact groups.}
The second step is to reduce the analysis of the induced
representation to a question about representations of compact
groups. To this end we fix a special maximal compact subgroup $G
\sub \bG$. From the Iwasawa decomposition $\bG = G\cdot \bP$ (cf.
\cite{Ca}, sec. 3.5), we deduce an isomorphism of
$G$-representations

$$\Ind^{\bG}_{\bP}(\rho)|_G
\stackrel{\simeq}{\lra} \Ind^G_{G \cap \bP}(\rho) \,.$$

\medskip

We let $I \sub G$ be an Iwahori subgroup, and put $P^+ =  G \cap \bP$. We note that every
element of the Weyl group $W$ of $\frG$ with respect to $\frS$ has
a representative in $G$ (cf. \cite{BT1}, 4.2.3, \cite{Ca}, p. 140 (b)),
and hence we
identify $W$ with $(N_\bG(\bS) \cap G)/(Z_\bG(\bS) \cap G)$. The
same remark applies to the Weyl group $W_\bP$. From the
Bruhat-Tits decomposition (cf. \cite{Ca} 3.5)  we deduce that

$${\bf G} = \coprod_{w \in W/W_\bP} Iw\bP \,.$$

It follows that

$$G = \coprod_{w \in W/W_{\bP}} IwP^+ \,.$$

\medskip

This in turn shows that, as representations of
$I$, we have

$$\Ind^G_{G \cap \bP}(\rho)|_I \cong
\bigoplus_{w \in  W/W_{\bP}}
\Ind^I_{P_w^+}(\rho^w) \,,$$

\medskip

where $P_w^+ = I \cap wP^+w^{-1}$. Here, we fix once and for all representatives for elements in $W/W_{\bP}$ as follows: first we take from each coset $wW_{\bP}$ its Kostant representatives, $\dot{w}$ say, and then we fix a representative of $\dot{w}$ in $G$ (that is possible, cf. above). When we write $w \in W/W_{\bP}$,  then $w$ denotes this specifically chosen representative in $G$. We let $M^w(\rho)$ be the $D(P,K)$-module dual to $\Ind^I_{P_w^+}(\rho^w)$, i.e.,

$$M^w(\rho) = \left(\Ind^I_{P_w^+}(\rho^w)\right)'_b
= D(I,K) \otimes_{D(P^+_w,K)} V'_w \,,$$

\medskip

where $V_w$ denotes the vector space $V$ equipped with the action $\rho^w$ of
$P^+_w$ given by $\rho^w(h) = \rho(w^{-1}hw)$, cf. (\ref{dual_iso}). We get

$$\left(\Ind^G_{G \cap \bP}(\rho)\right)'_b
= \bigoplus_{w \in  W/W_{\bP}}
M^w(\rho) \,.$$

\medskip

Assuming that all $D(I,K)$-modules $M^w(\rho)$ are simple we show
in Prop. \ref{non-isomorphic modules for different w} that they are
pairwise non-isomorphic. This implies that the left hand side is a
simple $D(G,K)$-module, which in turn shows that $\Ind^G_{G \cap
\bP}(\rho)$ is a topologically irreducible $G$-representation.
\end{para}

\bigskip

\begin{para}\label{induced rep parahoric}
{\it From $D(I,K)$-modules to $D_r(I,K)$-modules.}
Because the representation $V$ is finite-dimensional, the canonical map
$V_w' \ra D_r(P_w^+,K) \otimes_{D(P^+_w,K)} V_w'$ is an
isomorphism for $r<1$ sufficiently close to $1$, cf. Prop. \ref{extended action}.
Hence

$$\begin{array}{rl}
M^w_r(\rho) & := D_r(I,K) \, \otimes_{D(I,K)} M^w(\rho)
= D_r(I,K) \, \otimes_{D(P^+_w,K)} V'_w \\
 & \\
 & = D_r(I,K) \otimes_{D_r(P_w^+,K)}
\left( D_r(P_w^+,K) \otimes_{D(P^+_w,K)} V_w' \right) \\
 & \\
 & =  D_r(I,K) \otimes_{D_r(P_w^+,K)} V_w' \,.
\end{array}$$

\medskip

By \cite{ST2}, Lemma 3.9, $M^w(\rho)$ is a simple $D(I,K)$-module if
$M^w_r(\rho)$ is a simple $D_r(I,K)$-module for a sequence of
$r$'s tending to $1$, so that we are done if we show that all $M^w_r(\rho)$
are simple $D_r(I,K)$-modules, assuming the simplicity of $m(\rho)$
as an $U(\frg)$-module.

\medskip

In order to study the modules $M^w_r(\rho)$ we use the Iwahori product
decomposition, cf. Lemma \ref{Zerlegung},
$$I = U^-_w \cdot P^+_w \,,$$

\medskip

where $U^-_w = I \cap w\bU^-w^{-1}$ and $\bU^-$ is the unipotent
radical of the parabolic subgroup opposite to $\bP$.
We choose an open uniform normal subgroup $I_0 \sub I$ such that

\begin{equation*}
I_0 = U^-_{w,0} \cdot P^+_{w,0} \,\,\,\, \mbox{where} \,\,\,\,
U^-_{w,0} = U^-_w \cap I_0 \,\,\,\, \mbox{and} \,\,\,\, P^+_{w,0}
= P^+_w \cap I_0 \,.
\end{equation*}

\medskip

We will show that it is possible to choose $I_0$ in such a way that $P^+_{w,0}$
is uniform and $U^-_{w,0}$ is $L$-uniform (cf. Remark \ref{L_uniform}). The completed distribution algebras of $I$ ($U^-_w$ and $P^+_w$,
resp.) are defined by means of the canonical $p$-valuation on
$I_0$ ($U^-_{w,0}$ and $P^+_{w,0}$, resp.), cf.
\ref{norms}. For $r<1$ sufficiently close to $1$,
we have a canonical isomorphism of Banach spaces

$$D_r(U_w^-,K) \hat{\otimes}_K D_r(P_w^+,K) \stackrel{\simeq}{\lra}
D_r(I,K) \,,$$

\medskip

which in turn gives rise to a canonical isomorphism of
$D_r(U_w^-,K)$-modules

$$D_r(U_w^-,K) \otimes_K V_w' \stackrel{\simeq}{\lra} M^w_r(\rho) \,,$$

\medskip

cf. Prop. \ref{extended action}. Using the integrality of the
distribution algebra $D_r(U_w^-,K)$, cf. Prop. \ref{integrality of dist
algebra}, together with Corollary \ref{nontrivial intersection},
we prove that any non-zero $D_r(I,K)$-submodule $N$ of
$M^w_r(\rho)$ has non-zero intersection with

$$m^w_r(\rho) := U_r(\fru^-_w,U^-_{w,0}) \otimes_K V_w' \,,$$

\medskip

cf. Prop. \ref{from_Mr to mr}. Using general results about orthogonal
bases we can even infer that $N$ has non-zero intersection with

$$m^w(\rho) = U(\fru^-_w) \otimes_K V_w' \simeq U(\frg)
\otimes_{U(\frp^+_w)} V_w' \,.$$

\medskip

But it is not difficult to see that if $m(\rho) := m^1(\rho) =
U(\frg) \otimes_\frp V'$ is a simple $U(\frg)$-module, then
$m^w(\rho)$ is a simple $U(\frg)$-module for all $w$, and this
implies that $M^w_r(\rho)$ is a simple $D_r(I,K)$-module for all
$w$ and all $r$ sufficiently close to $1$. Hence, by our previous
remark, $M^w(\rho)$ is a simple $D(I,K)$-module for all $w$. From
what we have said at the end of \ref{Passage to representations of
compact groups} it then follows that $\Ind^G_{G \cap \bP}(\rho)$
is a topologically irreducible $G$-representation.
\end{para}

\bigskip

\subsection{Parahoric subgroups and their distribution algebras}
\label{Parahoric subgroups and their distribution algebras}

\begin{para}\label{defining parahoric subgroups}
For the following compare
\cite{Ca}, sec. 3.5. The torus $\frS$ determines an apartment
$\cA$ in the Bruhat-Tits building of $\frG$ over $L$. We fix a
special vertex $x_0$ in the apartment $\cA$. Then there is a
unique conical chamber $\cC$ in $\cA$ having $x_0$ as apex and
satisfying the following property: for every $u$ in the unipotent
radical $\bU_0$  of $\bP_0$ the intersection $\cC \cap u \cC$
contains a translate of $\cC$. Moreover, there is a unique chamber
$C_0$ in $\cC$ having $x_0$ as one of its vertices. We let $G
\subset \bG$ be the stabilizer of $x_0$ and $I \subset G$ the
pointwise stabilizer of $C_0$. Let $\bU^-$ be the
unipotent radical of the parabolic subgroup of $\bG$ opposite to
$\bP$. For $w \in  W / W_{\bP}$, we put

$$P_w^+ = I \cap w\bP w^{-1} \,, \,\, U_w^+ =  I \cap w\bU w^{-1} \,,
\,\, U_w^- = I \cap w\bU^-w^{-1} \,.$$
\end{para}




\bigskip

\begin{lemma}\label{Zerlegung} The multiplication map

$$P_w^+ \times U^-_w  \lra I$$

\medskip

is an isomorphism of locally $L$-analytic manifolds. In
particular, there are decompositions

$$I = P_w^+ \cdot U^-_w = U^-_w \cdot P_w^+ \,.$$

\medskip
\end{lemma}

{\it Proof.} Let $\bZ$ be the centralizer of $\bS$ in $\bG$ and denote by $Z_c$ the kernel of the natural homomorphism
$\bZ(L) \to X_\ast(S)_{\mathbb R}.$   For a (reduced) root $\alpha \in \Phi^{\red}$, let $\bU_\alpha$ be its root subgroup in $\bG$. Let $X_\alpha = \bU_\alpha(L) \cap I$.  By our assumption on $\frG$, the condition in \cite{Ti},  first sentence in 3.1.1, is met (cf. the sentence preceding sec. 3.1.1 in \cite{Ti}). Therefore, by the last sentence in \cite{Ti}, sec. 3.1.1, we have a product decomposition

$$I = Z_c \cdot\prod_{\alpha \in \Phi^{\red}} X_{\alpha}  \,,$$

where the product is taken with respect to any ordering on the
roots. The claim follows immediately.
\qed

\bigskip

\begin{para}\label{defining norms} {\it Defining norms on $D(I,K)$.}
We start with the construction of the subgroup $I_0 \sub I$
mentioned in \ref{induced rep parahoric}. It follows from
\cite{S-S}, Prop. I.2.7, that the subgroups $U_{x_0}^{(e)}$ of $G$ constructed there (we take for the facet $F$ in \cite{S-S} the special vertex $x_0$ of \ref{defining parahoric subgroups}) possess an Iwahori decomposition with respect to the parabolic subgroup $\bP$:

$$U_{x_0}^{(e)} = (U_{x_0}^{(e)} \cap \bU^-) \cdot (U_{x_0}^{(e)} \cap \bP) \,.$$

\medskip

But as the groups $U_{x_0}^{(e)}$ are normal in $G$ (\cite{S-S}, three lines before Prop. I.2.7),
they possess an Iwahori decomposition with respect to any parabolic subgroup of the form $w\bP w^{-1}$:

$$U_{x_0}^{(e)} = w(U_{x_0}^{(e)} \cap \bU^-)w^{-1} \cdot w(U_{x_0}^{(e)} \cap \bP)w^{-1} = (U_{x_0}^{(e)} \cap w\bU^- w^{-1}) \cdot (U_{x_0}^{(e)} \cap w\bP w^{-1})\,.$$

\medskip

The groups $U_{x_0}^{(e)}$ form a fundamental system of compact open neighborhoods of $1$ in $\bG$ (\cite{S-S}, Cor. I.2.9). Fix an $L$-uniform subgroup $H \sub I$ (cf. Remark \ref{L_uniform}). Then there is $e>0$ such that $U_{x_0}^{(e)} \sub H$. Let $\Lambda \sub \Lie(H) = \frg$ be the $\fro_L$ lattice from Lemma \ref{L_uniformgroups}, i.e., $\exp_H$ maps $\Lambda$ homeomorphically onto $H$. Define:

$$\begin{array}{lcl}
\Lambda_0 &:=& \exp_H^{-1}(U_{x_0}^{(e)})\\
&&\\
\Lambda_1 &:=& \exp_H^{-1}(U_{x_0}^{(e)} \cap w\bU^- w^{-1})\\
&&\\
\Lambda_2 &:=& \exp_H^{-1}(U_{x_0}^{(e)} \cap w\bP w^{-1})\\
\end{array}
$$

\medskip

Each $\Lambda_i$, $i=0,1,2$, is a $\Zp$-lattice in $\Lambda$. For $t \gg 0$ we then have that $p^t \Lambda_i$ is a uniform $\Zp$-Lie algebra in the sense of \cite{DDMS}, sec. 9.4. Fix such a $t$ and define

$$\begin{array}{lclclcl}
U^-_{w,0} &:=& \exp_H(p^t \Lambda_1) &\sub& U_{x_0}^{(e)} \cap w\bU^- w^{-1} &\sub& U^-_w \,,\\
&&&&&&\\
P^+_{w,0} &:=& \exp_H(p^t \Lambda_2) &\sub& U_{x_0}^{(e)} \cap w\bP w^{-1} &\sub& P^+_w \,,\\
&&&&&&\\
I_0 &:=& \exp_H(p^t \Lambda_0) &\sub& U_{x_0}^{(e)} &\sub& I \,.
\end{array}$$

\medskip

The subgroup $U^-_{w,0}$ ($P^+_{w,0}$, $I_0$, resp.) consists of the $p^t$-th powers of elements of $U_{x_0}^{(e)} \cap w\bU^- w^{-1}$ ($U_{x_0}^{(e)} \cap w\bP w^{-1}$, $U_{x_0}^{(e)}$, resp.), and is thus a characteristic subgroup of $U_{x_0}^{(e)} \cap w\bU^- w^{-1}$ ($U_{x_0}^{(e)} \cap w\bP w^{-1}$, $U_{x_0}^{(e)}$, resp.). As $U_{x_0}^{(e)}$ is normal in $G$, it follows that the group $U_{x_0}^{(e)} \cap w\bU^- w^{-1}$ ($U_{x_0}^{(e)} \cap w\bP w^{-1}$, $U_{x_0}^{(e)}$, resp.) is normal in $U^-_w$ ($P^+_w$, $I$, resp.). And this implies that $U^-_{w,0}$ ($P^+_{w,0}$, $I_0$, resp.) is normal in $U^-_w$ ($P^+_w$, $I$, resp.). Each group $U^-_{w,0}$, $P^+_{w,0}$, and $I_0$, is uniform pro-$p$.

\medskip

We now show that $U^-_{w,0}$ is actually $L$-uniform. Because of the root decomposition of $U_{x_0}^{(e)}$ in \cite{S-S}, Prop. I.2.7, it suffices to consider $U^-_{w,0} \cap \bU_{\alpha}$ for some $\alpha \in (\Phi^- \setminus \Phi_\bP) \cap \Phi^{{\rm red}}$. Here, $\Phi_\bP$ is the root system of the Levi subgroup of $\frP$ which contains $\frS$, and $\bU_\alpha$ is the generalized root group with Lie algebra $\frg_\alpha \oplus \frg_{2\alpha}$.\footnote{Here $\frg_{2\alpha}=0$ if $2\alpha$ is not a root.} It follows from \cite{BT2}, last sentence in sec. 5.2.2, that $U_{x_0}^{(e)} \cap \bU_\alpha$ is equal to $\cU(\fro_L)$, where $\cU$ is a group scheme over $\fro_L$ which is, as a scheme, isomorphic to ${\rm Spec}\left({\rm Sym}_{\fro_L}(N) \right)$ with a free $\fro_L$-module $N$ of finite rank. Hence we see that $U^-_{w,0} \cap \bU_{\alpha} = \cU_t(\fro_L)$, where $\cU_t = {\rm Spec}\left({\rm Sym}_{\fro_L}(N_t)\right)$ with $N_t = p^{-t}N$. The group of $\fro_L$-value
 d points $\cU_t(\fro_L)$ is identified with the $\fro_L$-algebra homomorphisms ${\rm Sym}_{\fro_L}(N_t) \ra \fro_L$, which are uniquely determined by restriction to $N_t \sub {\rm Sym}_{\fro_L}(N_t)$. We may assume that the zero section of $\cU_t$ corresponds to the map which sends $N_t$ to zero. Let $(N_t) \sub {\rm Sym}_{\fro_L}(N_t)$ be the ideal generated by $N_t$.  Then the Lie algebra of the group scheme $\cU_t$ can be identified with the relative tangent space $\Hom_{\fro_L}((N_t)/(N_t)^2, \fro_L)$. It is easy to see that the exponential map

$$\exp_{\cU_t(\fro_L)}: {\rm Lie}(\cU_t) \otimes_{\fro_L} L \dashrightarrow \cU_t(\fro_L) \,,$$

\medskip

when restricted to a sufficiently small submodule of $\rm{Lie}(\cU_t)$, corresponds then to the restriction (to that submodule) of the composed map

$$\begin{array}{lcl}
{\rm Lie}(\cU_t) = \Hom_{\fro_L}((N_t)/(N_t)^2, \fro_L) & \lra & \Hom_{\fro_L}(N_t, \fro_L) \\
&&\\
& \lra & \Hom_{\fro_L-\rm{algebras}}({\rm Sym}_{\fro_L}(N_t), \fro_L) = \cU_t(\fro_L) \,.
\end{array}$$

\medskip

The first map is defined by restricting a map $(N_t)/(N_t)^2 \ra \fro_L$ to $N_t \sub (N_t)/(N_t)^2$. The maps in this sequence are all bijections. Because $\cU_t(\fro_L)$ is uniform pro-$p$, the exponential map maps a lattice in ${\rm Lie}(\cU_t)$ bijectively onto $\cU_t(\fro_L)$, and this lattice must therefore be $\rm{Lie}(\cU_t)$, which is an $\fro_L$-module. Hence we conclude that $U^-_{w,0} \cap \bU_{\alpha} = \cU_t(\fro_L)$ is $L$-uniform, and it follows that $U^-_{w,0}$ is $L$-uniform too.

\medskip

We have $I_0= U^-_{w,0} \cdot
P^+_{w,0}$, and for the canonical $p$-valuation on $I_0$ it
follows from the identity $P_{i+1}(I_0)=\exp(p^i \cdot p^t\Lambda_0)$, together with \cite{DDMS}, Lemma 4.14 (iv), that for $x \in p^t \Lambda_1$ and $y \in p^t\Lambda_2$ we have

$$\omega^{can}(\exp(x)\exp(y)) = \vep_p + \min\{a,b\}
= \min\{\omega^{can}(\exp(x)), \omega^{can}(\exp(y))\} \,,$$

\medskip

where $a,b$ are such that $x \in p^a\cdot (p^t \Lambda_1) \setminus p^{a+1} \cdot (p^t \Lambda_1)$ and
$y \in p^b \cdot (p^t\Lambda_2) \setminus p^{b+1}\cdot (p^t\Lambda_2)$ (with $a$ (resp. $b$) being
$\infty$ if $x = 0$ (resp. $y = 0$)). In order to define the norms
$q_r$ on the rings $D(U^-_w,K)$, $D(P^+_w,K)$ and $D(I,K)$, we
work with the uniform normal subgroups $U^-_{w,0}$, $P^+_{w,0}$
and $I_0$, and the canonical $p$-valuations on these groups,
following the recipe explained in \ref{norms}.
\end{para}

In the following Proposition, all completed tensor products are meant to be the completions of the ordinary tensor products with respect to the projective tensor product topology, cf. \cite{S}, sec. 17. We remark that for Fr\'echet spaces the projective and the inductive tensor product topology coincide, cf. \cite{S}, Prop. 17.6.

\bigskip

\begin{prop}\label{decomposition of distring}
(i) The decompositions $I = U^-_w \cdot P^+_w$, $I_0 = U^-_{w,0}
\cdot P^+_{w,0}$ induce isomorphisms of topological $K$-vector
spaces

$$\begin{array}{c}
C^{an}_\Qp(U^-_w,K) \hat{\otimes}_{K}  C^{an}_\Qp(P^+_w,K)
\stackrel{\simeq}{\lra} C^{an}_\Qp(I,K) \\
\\
C^{an}_\Qp(U^-_{w,0},K) \hat{\otimes}_{K}  C^{an}_\Qp(P^+_{w,0},K)
\stackrel{\simeq}{\lra} C^{an}_\Qp(I_0,K) \\
\\
C^{an}_L(U^-_w,K) \hat{\otimes}_{K}  C^{an}_L(P^+_w,K)
\stackrel{\simeq}{\lra} C^{an}_L(I,K) \\
\\
C^{an}_L(U^-_{w,0},K) \hat{\otimes}_{K} C^{an}_L(P^+_{w,0},K)
\stackrel{\simeq}{\lra} C^{an}_L(I_0,K) \,.
\end{array}$$

\medskip

(ii) We equip the rings $D_\Qp(I_0,K)$, $D_\Qp(U^-_{w,0},K)$ and
$D_\Qp(P^+_{w,0},K)$ with the norm $\| \cdot \|_r$, $\frac{1}{p} <
r < 1$, associated to the canonical $p$-valuation. The rings
$D_\Qp(I,K)$, $D_\Qp(U^-_w,K)$ and $D_\Qp(P^+_w,K)$ carry the
maximum norms $q_r$, and $D_L(I,K)$, \linebreak $D_L(U^-_w,K)$,
$D_L(P^+_w,K)$ as well as $D_L(I_0,K)$, $D_L(U^-_{w,0},K)$ and
$D_L(P^+_{w,0},K)$ are equipped with the quotient norms
$\bar{q}_r$. On the tensor products of these spaces we put the
usual induced norm. Then the isomorphisms in (i) induce isometries
of topological $K$-vector spaces

$$\begin{array}{c}
D_\Qp(I,K) \stackrel{\simeq}{\lra} D_\Qp(U^-_w,K)
\hat{\otimes}_{K} D_\Qp(P^+_w,K) \\
\\
D_\Qp(I_0,K) \stackrel{\simeq}{\lra} D_\Qp(U^-_{w,0},K)
\hat{\otimes}_{K} D_\Qp(P^+_{w,0},K)\\
\\
D_L(I,K) \stackrel{\simeq}{\lra} D_L(U^-_w,K)
\hat{\otimes}_{K} D_L(P^+_w,K) \\
\\
D_L(I_0,K) \stackrel{\simeq}{\lra} D_L(U^-_{w,0},K)
\hat{\otimes}_{K} D_L(P^+_{w,0},K) \,.
\end{array}$$

\medskip

(iii) The isometries in (ii) furnish isometries of the completions

$$\begin{array}{c}
D_r(I,K) \stackrel{\simeq}{\lra} D_r(U^-_w,K)
\hat{\otimes}_{K} D_r(P^+_w,K) \\
\\
D_r(I_0,K) \stackrel{\simeq}{\lra} D_r(U^-_{w,0},K)
\hat{\otimes}_{K} D_r(P^+_{w,0},K) \,.
\end{array}$$

\end{prop}

{\it Proof.} (i) This follows from \cite{ST3}, Lemma A.1 and Prop. A.2.

(ii) Let $h_1, \ldots, h_{d'}$ be a $p$-basis  of $U^-_{w,0}$ and
$h_{d'+1}, \ldots, h_d$ be a $p$-basis of $P^+_{w,0}$. For $u \in
U^-_{w,0}$, $p \in P^+_{w,0}$, we have by the discussion in
\ref{defining norms}

$$\omega^{can}(up) = \min\{\omega^{can}(u), \omega^{can}(p)\} \,,$$

\medskip

so that $h_1, \ldots, h_{d'}, h_{d'+1}, \ldots, h_d$ is a
$p$-basis of $I_0$. Elements of $D_\Qp(I_0,K)$ have a unique
expansion as series of the form

$$\sum\nolimits_{n \in \bbN_0^d}d_n{\bf b}^n \,,$$

\medskip

with ${\bf b}^n = (h_1-1)^{n_1} \cdot \ldots \cdot (h_d-1)^{n_d}$.
From this we deduce immediately that the canonical map

$$D_\Qp(I_0,K) \stackrel{\simeq}{\lra} D_\Qp(U^-_{w,0},K)
\hat{\otimes}_K D_\Qp(P^+_{w,0},K)$$

\medskip

is an isometry when equipped with the norm $\| \cdot \|_r$ on the
left hand side and with the induced norm on the tensor product.
Consider the canonical commutative diagram

$$\begin{array}{ccc}
D_\Qp(I_0,K) & \stackrel{\simeq}{\lra} & D_\Qp(U^-_{w,0},K)
\hat{\otimes}_K D_\Qp(P^+_{w,0},K)\\
\downarrow & & \downarrow \\
D_L(I_0,K) & \stackrel{\simeq}{\lra} & D_L(U^-_{w,0},K)
\hat{\otimes}_K D_L(P^+_{w,0},K)
\end{array}$$

\medskip

Let $\bar{q}_r$ denote the quotient norm on $D_L(U^-_{w,0},K)$ as
well as on $D_L(P^+_{w,0},K)$, and let $\bar{q}_r \otimes
\bar{q}_r$ be the induced norm on the tensor product. On the other
hand, let $\| \cdot \|_r \otimes \| \cdot \|_r$ be the norm on
$D_\Qp(U^-_{w,0},K) \hat{\otimes}_K D_\Qp(P^+_{w,0},K)$, and
denote by $\overline{\| \cdot \|_r \otimes \| \cdot \|_r}$ the
norm on $D_L(U^-_{w,0},K) \hat{\otimes}_K D_L(P^+_{w,0},K)$
induced by the surjection

$$D_\Qp(U^-_{w,0},K)
\hat{\otimes}_K D_\Qp(P^+_{w,0},K) \twoheadrightarrow
D_L(U^-_{w,0},K) \hat{\otimes}_K D_L(P^+_{w,0},K) \,.$$

\medskip

By \cite{BGR}, Ch. 2.1, Prop. 6, we know that

$$\bar{q}_r \otimes \bar{q}_r =
\overline{\| \cdot \|_r \otimes \| \cdot \|_r} \,,$$

\medskip

which in turn shows that

$$D_L(I_0,K)  \stackrel{\simeq}{\lra}  D_L(U^-_{w,0},K)
\hat{\otimes}_K D_L(P^+_{w,0},K)$$

\medskip

is an isometry. Finally, using Remark \ref{equivalent_definition}
we can conclude that

$$D_L(I,K) \stackrel{\simeq}{\lra} D_L(U^-_w,K) \hat{\otimes}_K
D_L(P^+_w,K)$$

\medskip

is an isometry too.

(iii) This statement follows from (ii). \qed

\bigskip

\begin{prop}\label{integrality of dist algebra} (i) Let $\frH$ be a connected reductive group over $L$, $\frP \sub \frH$ a parabolic subgroup (defined over $L$), and let $\frN$ be the unipotent radical of $\frP$. Consider a compact open subgroup $H  \sub \frN(L)$. Choose an $L$-uniform subgroup $H_0 \sub H$ (cf. Remark \ref{L_uniform}), and define the norm $\bar{q}_r$ on $D(H,K)$ as explained in \ref{norms} (using the canonical $p$-valuation on $H_0$). Then there is a sequence $(r_m)_{m \ge 1}$ of numbers in $(\frac{1}{p},1) \cap p^\bbQ$, and converging to $1$, such that the ring $D_{r_m}(H,K)$ is an integral domain for all $m$.

\medskip

(ii) There is a sequence $(r_m)_{m \ge 1}$ of numbers in $(\frac{1}{p},1) \cap p^\bbQ$, and converging to $1$, such that the ring $D_{r_m}(U^-_w,K)$ is an
integral domain.
\end{prop}

{\it Proof.} (i) Let $\bar{q}_r$ be the norm on $D_r(H,K)$ defined
by means of the uniform subgroup $H_0$ (using the canonical
$p$-valuation on $H_0$) as explained in \ref{norms}. The key
idea of the proof is to embed $D_r(H,K)$ in another
distribution ring of the form $D_{r'}(H',K)$ which is an integral
domain. First we note that for any given compact open subgroup $H' \sub \frN(L)$ there is an element $s$ in a (maximal) torus $\frS(L)$ (which is chosen to normalize $\frN(L)$) such that $s H s^{-1}$ is contained in $H'$. We choose a compact open subgroup $H' \sub \frN(L)$ which is $L$-uniform (cf. Remark \ref{L_uniform}). $H$ and $sHs^{-1}$ being isomorphic locally $L$-analytic groups, it suffices to show the assertion for $sHs^{-1}$, i.e., we may assume without loss of generality that $H$ is already contained in a compact open subgroup $H'$ which is $L$-uniform. The embedding $H \hra H'$ gives rise to a continuous embedding of Fr\'echet spaces

$$D(H,K) \hra D(H',K) \,.$$

\medskip

Denote by $\bar{q}^{H'}_r$ the norms on $D(H',K)$ defined by means of the canonical $p$-valuation on $H'$ ($H'$ is uniform pro-$p$). The topology of the Fr\'echet spaces $D(H,K)$ and $D(H',K)$, is defined by the family of norms $(\bar{q}_r)_r$ and $(\bar{q}^{H'}_r)_r$, respectively. Therefore, for any
$r'$, there is an inequality of norms on $D(H,K)$

\begin{numequation}\label{ineq1}
\bar{q}_{r'}^{H'} \le c \bar{q}_r
\end{numequation}

\medskip

for $r$ sufficiently close to $1$. Next we consider the decomposition

$$D(H',K)  =  \bigoplus_{h \in H'/H} \delta_h \cdot D(H,K) \,,$$

\medskip

and denote by $\mu_r$ the maximum norm on $D(H',K)$ induced by the norm $\bar{q}_r$:

$$\mu_r\left(\sum_{h \in H'/H} \delta_h \cdot \lambda_h \right) = \max\{\bar{q}_r(\lambda_h) \midc h \in H'/H \} \,.$$

\medskip

The norms $\mu_r$ may not be algebra norms in general, but they define the Fr\'echet topology on $D(H',K)$, as well. Hence for a given $r_0$, there is an inequality

\begin{numequation}\label{ineq2}
\mu_{r_0} \le c' \bar{q}^{H'}_{r'}
\end{numequation}

\medskip

for $r'$ sufficiently close to $1$. Choose any $r_0 \in (p^{-\frac{1}{\vep_p (p-1)}},1)$. Then, by Prop. \ref{injectivity2}, there is $r' \in (\frac{1}{p},1) \cap p^\bbQ$ such that (\ref{ineq2}) holds, and such that $D_{r'}(H',K)$ is an integral domain. This shows that there is a continuous map

$$D_{r'}(H',K) \lra (D(H',K),\mu_{r_0})^{\wedge} \,,$$

\medskip

where the right hand side is the completion of $D(H,K)$ with respect to $\mu_{r_0}$. Next, choose $r \in (r_0, 1) \cap p^\bbQ$ such that (\ref{ineq1}) holds, and that, in addition, there is $m \in \bbN$ such that $s= r^{p^m}$ satisfies $s> \frac{1}{p}$ and $p^{-1/(p-1)-1/eq^{em}} \le s^{\vep_p} < p^{-1/(p-1)}$. Here, $e$ is the ramification index of $L$ over $\Qp$, and $q$ is the cardinality of the residue field of $L$. Note that we have $s < r_0 < r$ and hence a commutative diagram of continuous maps

$$\xymatrixcolsep{3pc}\xymatrix{D_r(H,K) \ar[r] \ar[d] & D_{r'}(H',K) \ar[d]\\
D_{r_0}(H,K) \ar[r] \ar[d] & (D(H',K),\mu_{r_0})^{\wedge}\\
D_s(H,K)  }$$

\medskip

The map $D_{r_0}(H,K) \ra (D(H',K),\mu_{r_0})^{\wedge}$ is clearly injective, because

$$(D(H',K),\mu_{r_0})^{\wedge} = \bigoplus_{h \in H'/H} \delta_h \cdot D_{r_0}(H,K) \,.$$

\medskip

Moreover, by definition of the norms $\bar{q}_r$ we have

\medskip

$$D_r(H,K) = \bigoplus_{h \in H/H_0} \delta_h \cdot D_r(H_0,K) \,.$$

\medskip

Because $r$ and $s$ are as in Prop. \ref{injectivity1}, the map $D_r(H_0,K) \ra D_s(H_0,K)$ is injective (cf. Prop. \ref{injectivity2}). Therefore, it follows from the remark just made that the map $D_r(H,K) \ra D_s(H,K)$ is injective as well. This implies that the map $D_r(H,K) \ra D_{r_0}(H,K)$ must be injective too. Hence the composite map (cf. the diagram above)

$$D_r(H,K) \lra D_{r_0}(H,K) \lra (D(H',K),\mu_{r_0})^{\wedge}$$

\medskip

is injective. Hence we see that the map $D_r(H,K) \lra D_{r'}(H',K)$ must be injective as well. We have chosen $r'$ such that $D_{r'}(H',K)$ is an integral domain. All maps under consideration are ring homomorphisms, and so $D_r(H,K)$ is an integral domain as well.

\medskip

(ii) We remark that in \ref{defining norms} we have defined $U^-_{w,0}$ in such a way that it is $L$-uniform. Hence the assertion follows immediately from (i). \qed

\bigskip

\begin{cor}\label{non-trivial intersection for U}
There is a sequence $(r_m)_m$ of numbers in $(\frac{1}{p},1) \cap p^\bbQ$, and converging to $1$, such that the following holds for any $m$: any non-zero left ideal $I$ of $D_{r_m}(U^-_w,K)$ has non-zero intersection with $U_{r_m}(\fru^-_w,U^-_{w,0})$.
\end{cor}

{\it Proof.} In \ref{defining norms} the
lattice $p^t\Lambda_1 \sub \fru^-_w$ has been defined so as to assure that $U^-_{w,0} = \exp(p^t\Lambda_1)$ is $L$-uniform (cf. Remark \ref{L_uniform}). Using Prop. \ref{integrality of dist algebra} we can apply Corollary
\ref{nontrivial intersection} whose assertion is exactly the claim
made above. \qed

\bigskip

\subsection{Modules for the completed distribution algebras}

\begin{para}\label{set-up representation}
{\it The modules $M^w(\rho)$ and $m^w(\rho)$.} Let

$$\rho: \bP \lra \bM \lra \GL(V)$$

\medskip

be the fixed locally analytic representation of $\bP$ on $V$
(equipped with the unique Hausdorff locally convex topology) from
\ref{induced representation}. The representation of
$P^+_w$ which we get by restriction and conjugation
will be denoted by $(\rho^w,V_w)$, cf. \ref{induced rep
parahoric}. Consider its (strong) dual $V_w' := (V_w)'_b$, which is a $D(P^+_w,K)$-module. We recall our convention to denote  Lie algebras of locally analytic groups  by their corresponding  gothic letters:

$$\frp^+_w = \Lie(P^+_w) = w \Lie(P^+) w^{-1} \,, \,\,
\fru^-_w = \Lie(U^-_w) = w \Lie(\bU^-) w^{-1} \,, \,\, \frs =
\Lie(\bS) \,.$$

\medskip

Via the embedding $U(\frp^+_w) \hra D(P^+_w,K)$ we view $V'_w$ as
a Lie algebra representation of $\frp^+_w$. We put

$$m^w(\rho) := U(\frg) \otimes_{U(\frp^+_w)}
V'_w \,\, \cong \,\, U(\fru^-_w)_K \otimes_K V'_w \,.$$

\medskip

Next, consider the induced locally analytic representation
$\Ind^I_{P^+_w}(\rho^w)$ and let

$$M^w(\rho) = (\Ind^I_{P^+_w}(\rho^w))'_b = D(I,K) \otimes_{D(P^+_w,K)} V'_w$$

\medskip

be the corresponding $D(I,K)$-module. By
Prop. \ref{decomposition of distring}, there is a canonical
isomorphism of $D(U^-,K)$-modules

$$M^w(\rho) = D(I,K) \otimes_{D(P^+_w,K)} V'_w \,\,
\simeq \,\, D(U^-_w,K) \otimes_K V'_w \,.$$

\medskip

Thus we see that the natural map $m^w(\rho) \ra M^w(\rho)$ is
injective.
\end{para}

\bigskip

\begin{prop}\label{extended action}
(i) For $r \in (\frac{1}{p},1)$ sufficiently close to $1$
the canonical map

$$V'_w \lra D_r(P^+_w,K) \otimes_{D(P^+_w,K)} V'_w \,\, ,
\,\, v \mapsto 1 \otimes v \,,$$

\medskip

is an isomorphism of $K$-vector spaces. Consequently, via
this isomorphism, we can extend the continuous operation of
$D(P^+_w,K)$ on $V'_w$ to a continuous operation of
$D_r(P^+_w,K)$ on $V'_w$. \\

(ii) When the finitely generated $D_r(I,K)$-module $D_r(I,K)
\otimes_{D_r(P^+_w,K)} V'_w$ is equipped with its natural
topology, the canonical map

$$D_r(U^-_w,K) \otimes_K V'_w \lra
D_r(I,K) \otimes_{D_r(P^+_w,K)} V'_w$$

\medskip

is a topological isomorphism.
\end{prop}

{\it Proof.} (i)  The $K$-vector space $V'_w$ is finite-dimensional
and hence finitely generated as $D(P^+_w,K)$-module. Thus $V_w$ is
a strongly admissible representation of $P^+_w$ in the sense of \cite{ST1},
\S 3. This implies that $V'_w$ is a co-admissible $D(P^+_w,K)$-module
in the sense of \cite{ST2}, \S 3 (cf.
\cite{ST2}, paragraph following Prop. 6.4). Hence

$$\phi: V'_w \lra \lim_{\stackrel{\longleftarrow}{r<1}}
\left( D_r(P^+_w,K) \otimes_{D(P^+_w,K)} V'_w \right) $$

\medskip

is an isomorphism of $K$-vector spaces. By \cite{ST2}, Theorem A, p. 152,
for every $r \in (\frac{1}{p},1)$, the canonical map

$$\phi_r: V'_w \lra D_r(P^+_w,K) \otimes_{D(P^+_w,K)} V'_w$$

\medskip

has dense image, and is therefore surjective. For $r_1 < r_2 <1$ the map
$\phi_{r_1}$ factors through $\phi_{r_2}$. This implies that the family
$\ker(\phi_r)$ is decreasing with increasing $r$. The intersection of all
$\ker(\phi_r)$ is the zero space because $\phi$ is injective. As $V'_w$ is
finite-dimensional, there must be some $r \in (\frac{1}{p},1)$ such that
$\ker(\phi_r)$ is zero. Hence, there is necessarily some
$r \in (\frac{1}{p},1)$ such that
the map $V'_w \lra D_r(P^+_w,K) \otimes_{D(P^+_w,K)} V'_w$ is bijective.

\medskip

(ii) By Prop. \ref{decomposition of distring} it follows that the
map in (ii) is a bijection. If we give the right hand side its
natural $D_r(I,K)$-module topology, i.e., the quotient topology
induced by an arbitrary surjection of $D_r(I,K)$-modules

$$D_r(I,K)^m \rightarrow D_r(I,K) \otimes_{D_r(P^+_w,K)} V'_w
\,,$$

\medskip

then the map in (ii) is a continuous bijective map, and hence, by
the open mapping theorem, a homeomorphism. \qed

\bigskip

\begin{para}\label{M_r(rho)}
{\it The modules $M^w_r(\rho)$ and $m^w_r(\rho)$.}
For the rest of this section, $r<1$ denotes a real number
sufficiently close to 1 such that the assertions of \ref{extended action}
and \ref{non-trivial intersection for U} do hold. We put

$$M^w_r(\rho) := D_r(I,K) \otimes_{D(I,K)} M^w(\rho) \,,$$

\medskip

which we consider as a module over the
Banach algebra $D_r(I,K)$. By the foregoing proposition, we have

$$\begin{array}{rl}
D_r(I,K) \, \otimes_{D(I,K)} M^w(\rho) &
= D_r(I,K) \, \otimes_{D(P^+_w,K)} V'_w \\
 & \\
 & = D_r(I,K) \otimes_{D_r(P_w^+,K)}
\left( D_r(P_w^+,K) \otimes_{D(P^+_w,K)} V_w' \right) \\
 & \\
 & =  D_r(I,K) \otimes_{D_r(P_w^+,K)} V_w' \,.
\end{array}$$

\medskip

Using this and Prop. \ref{decomposition of distring} we see that the map

$$D_r(U^-_w,K) \otimes_K V_w' \lra M^w_r(\rho)$$

\medskip

is an isomorphism of $D_r(U^-_w,K)$-modules. Therefore, the natural map

$$M^w(\rho) \lra  M^w_r(\rho)$$

\medskip

is an embedding. We denote the topological
closure of $m^w(\rho)$ in $M^w_r(\rho)$ by $m^w_r(\rho)$. Because
$U_r(\fru^-_w,U^-_{w,0})$ is the topological closure of
$U(\fru^-_w)_K$ in $D_r(U^-_w,K)$ we have

$$m^w_r(\rho) = U_r(\fru^-_w,U^-_{w,0}) \otimes_K V'_w \,.$$

\medskip

Since the maximal split torus $\bS$ is contained in $\bP$, we have
a natural diagonal action of $U(\frs)$ on $m^w(\rho)$.
\end{para}

\medskip
\begin{lemma}\label{finite-dim'l weight spaces} (i) For any $\frx \in \frs$ the action of $\frx$ on $m^w(\rho)$ is continuous with respect to the topology on $m^w_r(\rho)$, and hence extends continuously to $m^w_r(\rho)$.

\medskip

(ii) The $U(\frs)$-module structure on $m^w(\rho)$ extends to a $U(\frs)$-module structure on $m^w_r(\rho)$.

\medskip

(iii) For any $\lambda \in \frs^\ast$ the weight space
$m^w_r(\rho)_\lambda$ is finite-dimensional.
\end{lemma}

{\it Proof.} We begin by considering again the way we have defined the subgroup
$U^-_{w,0}$ in \ref{defining norms}. Recall that we started with an $L$-uniform (normal) subgroup $H$ of $G$, and a lattice $\Lambda \sub \Lie(H)$ on which $\exp_H$ is defined and maps $\Lambda$ bijectively to $H$. Then we have chosen $e$ large enough such that $U^{(e)}_{x_0}$ is contained in $H$. Then we have put

$$\Lambda_1 := \exp_H^{-1}(U^{(e)}_{x_0} \cap w\bU^- w^{-1}) = \mbox{Ad}(w)\left(\exp_H^{-1}(U^{(e)}_{x_0} \cap \bU^-) \right)$$

\medskip

and $U^-_{w,0} = \exp_H(p^t\Lambda_1)$ for $t \gg 0$.

\bigskip

Our aim is to show that $\Lambda_1$ has an $\fro_L$-basis consisting of root elements.

\bigskip

By \cite{S-S}, Prop. I.2.7, we have

$$U^{(e)}_{x_0} \cap \bU^- = \prod_{\alpha \in (\Phi^- \setminus \Phi_\bP) \cap \Phi^{\rm{red}}} \left(U_{f^*_{x_0} + e} \cap U_\alpha\right)$$

\medskip

where $\Phi_\bP$ is the root system of the Levi subgroup of $\bP$ which contains the torus $\bS$, and where

\begin{numequation}\label{deco1}
U_{f^*_{x_0} + e} \cap U_\alpha = U_{\alpha, f^*_{x_0}(\alpha)+e} \,\, . \, U_{2\alpha,2f^*_{x_0}(\alpha)+e} \,.
\end{numequation}

\medskip

We refer to \cite{S-S} (and \cite{BT2}) for the definitions of the various groups appearing here. We want to show that

\begin{numequation}\label{deco2}
\begin{array}{l}
\exp_H^{-1}\left(U_{f^*_{x_0} + e} \cap U_\alpha\right) \\
\\
= \left(\exp_H^{-1}\left(U_{f^*_{x_0} + e} \cap U_\alpha\right) \cap \frg_\alpha \right) \oplus \left(\exp_H^{-1}\left(U_{f^*_{x_0} + e} \cap U_\alpha\right) \cap \frg_{2\alpha}\right) \,.
\end{array}
\end{numequation}

\medskip

Here $\frg_\alpha$, $\frg_{2\alpha}$ denote, as usual, the root subspaces.\footnote{Note that, according to the conventions used in \cite{S-S} and \cite{BT2}, one has $\Lie(U_\alpha) = \frg_\alpha \oplus \frg_{2\alpha}$.}
By identity (\ref{deco1}) it follows that (\ref{deco2}) holds if we have

\begin{numequation}\label{deco3}
\exp_H^{-1}\left(U_{\alpha, f^*_{x_0}(\alpha)+e}\right) = \left(\exp_H^{-1}\left(U_{\alpha, f^*_{x_0}(\alpha)+e}\right) \cap \frg_\alpha \right) \oplus \left(\exp_H^{-1}\left(U_{\alpha, f^*_{x_0}(\alpha)+e}\right) \cap \frg_{2\alpha} \right) \,.
\end{numequation}

\medskip

To show that (\ref{deco3}) holds is a non-trivial problem only if $2\alpha$ is a root, which we assume from now on. At this point we have to consider the actual definition of the groups involved, and this means to trace through the {\oe}uvre of Bruhat and Tits. We will give the reader only some guidelines where to look in \cite{BT2}. By \cite{BT2}, 5.1.16, one can easily reduce the problem to the case when $\frG$ is quasi-split. By \cite{BT2}, 4.3.5, one has $U_{\alpha,k} = \frU_{\alpha,k}(\fro_L)$, and this group has an explicit description in terms of a group denoted $H(L,L_2)$ in \cite{BT2}, cf. \cite{BT2}, 4.1.9, 4.1.15.\footnote{The field denoted $L$ in \cite{BT2}, ch. 4, is not the same as our field $L$; our field $L$ is the field $K$ in \cite{BT2}, ch. 4.} But the description of $H(L,L_2)$ given in \cite{BT2}, 4.1.15, shows that (\ref{deco3}) does indeed hold.

\medskip

Hence we see that $\exp_H^{-1}(U^{(e)}_{x_0} \cap \bU^-)$ has an $\fro_L$-basis consisting of weight vectors, and the same is therefore true for the lattice $\Lambda_1$ and for $p^t\Lambda_1$. Thus, let $(\frX_1, \ldots, \frX_d)$ be an $\fro_L$-basis of $p^t\Lambda_1 = \exp_H^{-1}(U^-_{w,0})$ which consists of weight vectors $\frX_i \in \frg_{\alpha_i}$.

\medskip

(i) Because $U^-_{w,0}$ is $L$-uniform, we can apply part (ii) of Kohlhaase's theorem \ref{Kohlhaase} to the basis $(\frX_1, \ldots, \frX_d)$: the closure $U_r(\fru^-_w,U^-_{w,0})$ of $U(\fru^-_w)
\otimes_L K$ in $D_r(U^-_{w,0},K)$, cf. \ref{notation closure},
consists exactly of those series

$$\sum\nolimits_{n \in \bbN_0^d} d_n \frX^n$$

\smallskip

for which

$$\lim\nolimits_{|n| \ra \infty} |d_n| \nu_r(\frX^n) = 0 \,.$$

\medskip

Here, $\nu_r$ is a norm on $U_r(\fru^-_w,U^-_{w,0})$ which is
equivalent to $\bar{q}_r$ and $\nu_r(\sum_n d_n \frX^n) = \sup_n
|d_n| \nu_r(\frX^n)$. Let $v_1, \ldots, v_k$ be a basis of $V'_w$
which consists of weight vectors: $\frx \cdot v_j = \gamma_j(\frx)
v_j$ for all $\frx \in \frs$, $1 \le j \le k$. Then we have for any $\frx \in \frs$:

$$\frx \cdot (\frX^n \otimes v_j) = \left(\gamma_j + \sum\nolimits_{1 \le i \le k} n_i \alpha_i\right)(\frx) \cdot (\frX^n \otimes v_j) \,.$$

\medskip

It is obvious that there is $C>0$ (depending on $\frx$) such that for all $j \in \{1, \ldots, k\}$ and all $n \in \bbN_0^d$ one has

$$\Big|\left(\gamma_j + \sum\nolimits_{1 \le i \le k} n_i \alpha_i\right)(\frx)\Big|_K \le C \,.$$

\medskip

This shows that the action of $\frx$ on $m^w(\rho)$ extends to a continuous endomorphism of $m^w_r(\rho)$.

\medskip

(ii) Is an immediate consequence of (i).

\medskip

(iii) In (i) we have seen that the elements $\frX^n \otimes v_j \in m^w_r(\rho)$ are weight vectors for the action of $\frs$. If

$$\mu = \sum\nolimits_{n \in \bbN_0^d, 1 \le j \le k}
d_{n,j} \frX^n \otimes v_j\; \in \; m^w_r(\rho) = U_r(U^-_0,K)
\otimes_K V'_w$$

\medskip

is an element of $m^w_r(\rho)_\lambda$, then

$$\lambda = \gamma_j + \sum\nolimits_{1 \le i \le k} n_i \alpha_i$$

\medskip

if $d_{n,j} \neq 0$. This is an immediate consequence of the fact
that the monomials $\{\frX^n\}_{n \in \bbN_0^d}$ are an orthogonal
basis of $U_{r'}(U^-_0,K)$  with respect to the norm $\nu_r$, cf.
Thm. \ref{Kohlhaase}. Because the characters $\alpha_i$ do all occur in
the Lie algebra $\fru^-_w$, there can only be finitely many
possibilities to write $\lambda$ as a sum as above. Therefore, the
weight spaces $m^w_r(\rho)$ are all finite-dimensional. \qed

\bigskip

{\it From $m^w_r(\rho)$ back to $m^w(\rho)$.} From the result
above we can deduce that closed $U(\frs)$-invariant subspaces of
$m^w_r(\rho)$ are in bijection with $U(\frs)$-invariant subspaces
of $m^w(\rho)$.

\bigskip

\begin{prop}\label{non-trivial intersection with U(s)-invariant
submodule} (i) We have an inclusion preserving bijection

\begin{eqnarray*}
\Big\{\begin{array}{c}
\mbox{closed $U(\frs)$-invariant} \\
\mbox{subspaces of $m^w_r(\rho)$}
\end{array} \Big\} &
\stackrel{\sim}{\longrightarrow} &
\Big\{ \mbox{$U(\frs)$-invariant subspaces of $m^w(\rho)$} \Big\} \\
W & \longmapsto & W\cap m^w(\rho)
\end{eqnarray*}

\medskip

Any $U(\frs)$-invariant subspace $W$ of $m^w_r(\rho)$ is the
direct sum of its weight components: $W = \oplus W_\lambda$.

\medskip

(ii) Let $N \subset m^w_r(\rho)$ be a closed $U(\frs)$-invariant
subspace, $N \neq 0$. Then $N \cap m^w(\rho) \neq 0$. In
particular, any weight vector for the action of $\frs$ lies
already in $m^w(\rho)$.
\end{prop}

{\it Proof.} (i) This statement follows from Lemma
\ref{finite-dim'l weight spaces} and \cite{Fe} 1.3.12.

(ii) This is an immediate consequence of (i). \qed

\bigskip

\begin{lemma}\label{density} There is $r' \in (0,1)$ such that
$U(\fru^-_w) \otimes_L K$ is dense in $D_{r'}(U^-_{w,0},K)$.
\end{lemma}

{\it Proof.} Let $(\frx_i)_i$ be a basis of the $\bbZ_p$-Lie algebra $p^t\Lambda_1$
of $U^-_{w,0}$, with the notation as introduced in sec.
\ref{defining norms}. Then the elements $(\exp(\frx_i))_i$,
for an arbitrary but fixed ordering, form a $p$-basis of $U^-_{w,0}$.
Put $b_i = \exp(\frx_i)-1 \in D(R^L_\Qp U^-_{w,0},K)$,
so that the monomials ${\bf b}^n$ in the $b_i$ form an orthogonal basis
of $D_r(R^L_\Qp U^-_{w,0},K)$, as in sec. \ref{norms}.
For $0<r'<1$ sufficiently close to $0$ the right hand side of the identity
$b_i = \exp(\frx_i)-1$ can be expanded as a convergent series

$$\sum\nolimits_{\nu \ge 1} \frac{\frx_i^\nu}{\nu !}$$

\medskip

in $D_r(R^L_\Qp U^-_{w,0},K)$. This shows that we can approximate
$b_i$ by elements in $U(\fru^-_w)_K$. Hence $U(\fru^-_w) \otimes_L K$
is dense in $D_{r'}(U^-_{w,0},K)$. \qed

\bigskip

\begin{prop}\label{from_Mr to mr} Every non-trivial $U(\frs)$-invariant
$D_r(U^-_w,K)$-submodule of \linebreak $M^w_r(\rho)$ has a non-trivial
intersection with $m^w_r(\rho)$.
\end{prop}

{\it Proof.} (See also Proposition 11 of \cite{Fr}.) Fix a basis
$v_1, \ldots, v_k$ of weight vectors in $V'_w$ with respect to the
action of $U(\frs)$. We obtain $D_r(U^-_w,K)$-submodules

$$D_r(U^-_w,K) \otimes Kv_i \sub M^w_r(\rho)$$

\medskip

which are $U(\frs)$-invariant. Consider the projections

$$\pr_i: M^w_r(\rho) = D_r(U^-_w,K) \otimes_K V'_w
\ra D_r(U^-_w,K) \otimes_K Kv_i$$

\medskip

and let $N \sub M^w_r(\rho)$ be a $D_r(U^-_w,K)$ submodule which
is $U(\frs)$-invariant. By defining

$$N^{(i)}= \bigcap_{1\leq j <i} \ker(\pr_j) \cap N \,,$$

\medskip

we obtain a descending filtration of $D_r(U^-_w,K) \times
U(\frs)$-modules

$$0 = N^{(k+1)} \sub N^{(k)} \sub \ldots \sub N^{(2)}
\sub N^{(1)} = N \,.$$

\medskip

Let $1 \le i \le k$ be the unique index with

$$0 = N^{(i+1)} \sub N^{(i)} \neq 0 \,.$$

\medskip

Identifying $D_r(U^-_w,K) \otimes_K Kv_i$ with $D_r(U^-_w,K)$ as
$D_r(U^-_w,K)$-modules we see by Corollary \ref{non-trivial
intersection for U} that

$$\pr_i(N^{(i)}) \cap (U_r(\fru^-_w,U^-_{w,0}) \otimes_K Kv_i) \neq 0.$$

\medskip

By applying Prop. \ref{non-trivial intersection with
U(s)-invariant submodule} we can infer that

$$\pr_i(N^{(i)}) \cap (U(\fru^-_w) \otimes_L Kv_i) \neq 0 \,.$$

\medskip

Therefore, there is an element $F \in N^{(i)}$ such that $\pr_i(F)
\in U(\fru^-_w) \otimes_L Kv_i$ is a weight vector (again by
Prop. \ref{non-trivial intersection with U(s)-invariant submodule}). In
order to prove the statement of our proposition, it suffices to
show that $\pr_j(F) \in U_r(\fru^-_w,U^-_{w,0}) \otimes_K Kv_j$, $j
= i+1, \ldots, k$. Suppose that there is an index $i < j \leq k$
such that $\pr_j(F) \not\in U_r(\fru^-_w,U^-_{w,0}) \otimes_K
Kv_j$.

\medskip

We want to show that $\pr_j(F)$ can not be a weight vector for the
action of $U(\frs)$.

\medskip

To this end, choose $r' \in (0,r)$ sufficiently small such that
$U(\fru^-_w) \otimes_L K$ is dense in
$D_{r'}(U^-_{w,0},K)$, cf. Lemma \ref{density}. Then we have
a commutative diagram of embeddings

$$\begin{array}{ccc}
U_{r'}(\fru^-_w,U^-_{w,0}) \otimes_K Kv_j & \stackrel{=}{\lra} &
D_{r'}(U^-_{w,0},K) \otimes_K Kv_j \\
\cup & & \cup \\
U_r(\fru^-_w,U^-_{w,0}) \otimes_K Kv_j & \hra & D_r(U^-_{w,0},K)
\otimes_K Kv_j
\end{array}$$

\medskip

Therefore, we can consider $\pr_j(F)$ as an element of
$U_{r'}(\fru^-_w,U^-_{w,0}) \otimes_K Kv_j$. If $\pr_j(F)$ was a
weight vector, it would then automatically be an element of
$U(\fru^-_w) \otimes_L Kv_j$, by Prop. \ref{non-trivial intersection
with U(s)-invariant submodule} (ii). Hence, a fortiori, $\pr_j(F)$
would be in $U_r(\fru^-_w,U^-_{w,0}) \otimes_K Kv_j$.

\medskip

Thus we have shown that $\pr_j(F)$ is not a weight vector. Hence we
may choose $\lambda \in U(\frs)$ such that $\lambda \cdot  \pr_j(F)$
is not a scalar multiple of $\pr_j(F)$. Let $C_\lambda \in K$ be
the scalar with

$$\lambda \cdot \pr_i(F) = C_\lambda \cdot \pr_i(F) \,.$$

\medskip

Then the non-zero element $(\lambda - C_\lambda) \cdot F$ is
contained in $N^{(i+1)} = 0$, which is a contradiction. \qed

\bigskip

\begin{cor}\label{simple lie algebra module implies simple M-r-module}
If $m(\rho) = U(\frg) \otimes_{U(\frp^+)} V'$ is a simple
$U(\frg)$-module then:

\medskip

(i) $m^w(\rho) = U(\frg) \otimes_{U(\frp^+_w)} V'_w$ is a simple
$U(\frg)$-module for every $w$.

\medskip

(ii) $M^w_r(\rho)$ is a simple $D_r(I,K)$-module for every $w$.
\end{cor}

{\it Proof.} (i) Note that the map

$$U(\frg) \otimes_{U(\frp^+)} V' \lra U(\frg) \otimes_{U(\frp^+_w)}
V'_w \,, \,\, \frz \otimes v \mapsto \mbox{ad}(w)(\frz) \otimes v
\,,$$

\medskip

is an isomorphism of the underlying vector spaces. It sends
$U(\frg)$-submodules to $U(\frg)$-submodules. The left hand side
is therefore a simple $U(\frg)$-module if and only if the right
hand side is a simple $U(\frg)$-module.

\medskip

(ii) Let $N \sub M^w_r(\rho)$ be a $D_r(I,K)$-submodule. It is
automatically closed because $D_r(I,K)$ is noetherian and
$M^w_r(\rho)$ is a finitely generated $D_r(I,K)$-module. By the
previous proposition, we see that $N \cap\, m^w_r(\rho) \neq 0$. By
Prop. \ref{non-trivial intersection with U(s)-invariant submodule} we
get that $N \cap\, m^w(\rho) \neq 0$. Since $m^w(\rho)$ is a simple
$U(\frg)$-module, we obtain an inclusion $m^w(\rho) \sub\, N$. Thus
we conclude that $1 \otimes V'_w$ is contained in $N$ and
therefore $N = M^w_r(\rho)$. \qed

\bigskip

\begin{thm}\label{simplicity for parahoric}
If $m(\rho) = U(\frg) \otimes_{U(\frp^+)} V'$ is a simple
$U(\frg)$-module then

$$M^w(\rho) = \left(\Ind^I_{P_w^+}(\rho^w)\right)'_b$$

\medskip

is a simple $D(I,K)$-module for every $w$.
\end{thm}

{\it Proof.} It follows from the definition of $M^w(\rho)$
that this is a finitely generated $D(I,K)$-module (using that $V$ is
finite-dimensional over $K$). Therefore, $\Ind^I_{P_w^+}(\rho^w)$
is a strongly admissible representation of $I$ in the sense of \cite{ST1},
\S 3. This implies that $V'_w$ is a co-admissible $D(P^+_w,K)$-module
in the sense of \cite{ST2}, \S 3 (cf.
\cite{ST2}, paragraph following Prop. 6.4). Now we can use \cite{ST2},
Lemma 3.9. and Cor. \ref{simple lie algebra module implies simple M-r-module}
to conclude that $M^w(\rho)$ is a simple $D(I,K)$-module. \qed

\bigskip

\subsection{The main result}  The last essential step, to show the
topological irreducibility of the induced representation

$$\Ind^{\bG}_{\bP}(\rho) \,\, \cong \,\, \Ind^G_{P^+}(\rho)$$

\medskip

is to prove that the various $D(I,K)$-modules $M^w(\rho)$ are
pairwise non-isomorphic.

\medskip

\begin{prop}\label{non-isomorphic modules for different w}
Every homomorphism $M^{w'}(\rho) \ra M^w(\rho)$ of
$D(I,K)$-modules for $w \neq w'$ is zero.
\end{prop}

{\it Proof.} Our proof is a slight generalization of \cite{Fr},
Proposition 12. A homomorphism $M^{w'}(\rho) \ra M^w(\rho)$
corresponds by duality, cf. \ref{Locally analytic
representations}, to a homomorphism

$$\Ind^I_{P^+_w}(\rho^w) \lra \Ind^I_{P^+_{w'}}(\rho^{w'})$$

\medskip

of locally analytic $I$-representations. By Frobenius reciprocity
this corresponds to a continuous $P^+_{w'}$-homomorphism

$$\Ind^I_{P^+_w}(\rho^w) \rightarrow V_{w'}=V \,.$$

\medskip

From the decomposition $I = U^-_w \cdot P^+_w$ we deduce an
isomorphism

$$\Ind^I_{P^+_w}(\rho^w) \stackrel{\simeq}{\lra}
C_L^{an}(U^-_w,V)$$

\medskip

of representations of $U^-_w \cap P^+_{w'}$. Since $w \neq w'$ the
intersection $w' \bU (w')^{-1} \cap w \bU^- w^{-1}$ contains a
root group. Here, $\bU$ is the unipotent radical of $\bP$.
Therefore,

$$U := w' \bU (w')^{-1} \cap U^-_w \cap P^+_{w'}$$

\medskip

is a non-trivial $L$-analytic group of positive dimension. The
group $w' \bU (w')^{-1}$ acts trivially on $V_{w'}$, by
assumption. Thus $U$ acts trivially on $V_{w'}$. A homomorphism

$$\Ind^I_{P^+_w}(\rho^w) \lra \Ind^I_{P^+_{w'}}(\rho^{w'})$$

\medskip

gives therefore rise to a continuous map

$$\phi: C_L^{an}(U^-_w,V) \lra V$$

\medskip

which is $U$-equivariant, with $U$ acting trivially on $V$. The
canonical projection $U^-_w \ra U \bksl U^-_w$ has a locally
$L$-analytic section, so that we can find an isomorphism of
locally $L$-analytic manifolds $U^-_w \simeq U \times U'$ with
some compact $L$-analytic manifold $U'$. This isomorphism we may
assume to be compatible with the action of $U$ by left translation
(acting trivially on $U'$). This in turn gives rise to an
isomorphism

$$C_L^{an}(U^-_w,V) \stackrel{\simeq}{\lra}
C_L^{an}(U,K) \hat{\otimes}_{K,\pi} C_L^{an}(U',V) \,,$$

\medskip

cf. \cite{ST3}, Lemma A.1 and Prop. A.2. Let $\lambda: V \ra K$ be a linear form and fix $g \in
C_L^{an}(U',V)$. The map

$$\phi_{\lambda,g}:  C_L^{an}(U,K) \ra K \,, \,\,
f \mapsto  \phi_{\lambda,g}(f) = \lambda(\phi(fg)) \,,$$

\medskip

then has the property that for all $u \in U$:

$$\phi_{\lambda,g}(x \mapsto f(ux)) = \phi_{\lambda,g}(f) \,.$$

\medskip

Because $U$ does not have a $p$-adic Haar measure,
there is no non-zero $U$-invariant continuous linear form
on $C_L^{an}(U,K).$ Hence, we find that $\phi_{\lambda,g} = 0$. The
functions of the form $f \cdot g$, $f \in C_L^{an}(U,K)$, $g \in
C_L^{an}(U',V)$, span a dense subspace of $C_L^{an}(U^-_w,V)$.
Hence we deduce that $\lambda \circ \phi = 0$ for all linear forms
$\lambda$ on $V$. This shows that $\phi$ is necessarily zero. \qd

\bigskip

\begin{thm}\label{main thm}
Suppose $m(\rho)$ is simple as a $U(\frg)$-module. Then

\medskip

(i) $M(\rho) = \left(\Ind^G_{P^+}(\rho)\right)'_b$ is simple as a
$D(G,K)$-module. A fortiori, the representation
$\Ind^G_{P^+}(\rho)$ is topologically irreducible.

\medskip

(ii) $\Ind^{\bG}_{\bP}(\rho)$ is a topologically irreducible
representation of $\bG = \frG(L)$.

\medskip

\end{thm}

{\it Proof.} Consider the decomposition

$$\left(\Ind^G_{P^+}(\rho)\right)'_b =
\bigoplus_{W/W_{\bP}} M^w(\rho) \,.$$

\medskip

Each term $M^w(\rho)$ is by Theorem \ref{simplicity for parahoric}
a simple $D(I,K)$-module. Furthermore, the summands are not
pairwise isomorphic by the previous proposition. Since the Weyl
group permutes the summands transitively, $M(\rho)$ is a simple
$D(G,K)$-module. By the relation between $D(G,K)$-modules and
representations of $G$, we conclude that $\Ind^G_{P^+}(\rho)$ is a
topologically irreducible representation of $G$, cf. \ref{Locally
analytic representations}. The second assertion follows
immediately from (i). \qed

\section{Examples}

\subsection{Irreducibility of Verma modules}

Theorem \ref{main thm} is of course only useful if one can determine
whether the generalized Verma module
$m(\rho) = U(\frg) \otimes_{U(\frp)} \rho'$ is irreducible or not.
To simplify this problem, one may pass to a finite extension $\tilde{K}$
of $K$ such that the Lie algebra

$$\frg_{\tilde{K}} = \frg \otimes_L \tilde{K}$$

\medskip

is split over $\tilde{K}$.
Clearly, if $m(\rho) \otimes_K \tilde{K}$ is irreducible as module
over $U(\frg_{\tilde{K}})$, then $m(\rho)$ is simple as $U(\frg)$-module.
If one can reduce the problem to a question about classical Verma modules,
i.e. those which are induced from one-dimensional representations of a
Borel subalgebra, one is in a particular simple situation. In the next
paragraph we recall the set-up from \cite{D}, Ch.7.

\medskip

Suppose $\frg$ is a split reductive Lie algebra over $K$.
Let $\frh$ be Cartan subalgebra of $\frg$,
$\frb$ be a Borel subalgebra containing $\frh$, $\Phi$ the root system of
$\frg$ with respect to $\frh$, and $\Phi^+ \supset \Delta$
the set of positive and simple roots, respectively. Denote by
$\delta = \frac{1}{2}\sum_{\alpha \in \Phi^+} \alpha$ the sum of positive
roots, as usual.
Let $\{X_\beta, \beta \in \Phi; H_\alpha, \alpha \in \Delta \}$
be a Chevalley basis of $\frg_{\der} = [\frg,\frg]$.
For a character $\lambda \in \Hom_K(\frh,K)$, let $K_\lambda$ be the
one-dimensional representation of $\frh$ on $K$ defined by $\lambda$
and extend this to a representation of $\frb$ by letting the nilpotent radical act
trivially.

\bigskip

\begin{thm}\label{BGG} (Bernstein, Gelfand, Gelfand)
The Verma module $U(\frg) \otimes_{U(\frb)} K_\lambda$ is a
simple $U(\frg)$-module if and only if $(\lambda+\delta)(H_\alpha) \notin
\bbZ_{>0}$ for all $\alpha \in \Delta$.
\end{thm}

We refer to \cite{D}, Thm. 7.6.24. Note that Dixmier uses the normalized
induction, i.e. his $\lambda$ is our $\lambda+\delta$.

\bigskip

\begin{example} Let $\frG = \GL_{2,L}$ be the group $\GL_2$ over $L$,
$\frP$ the Borel subgroup
of upper triangular matrices, $\frS \sub \frP$ the diagonal torus and

$$\chi: \bS = \frS(L) \ra K^* \,\,, \,\,
\left(\begin{array}{cc}
t_1 & 0 \\
0   & t_2
\end{array} \right) \mapsto \chi_1(t_1) \chi_2(t_2) \,,$$

\medskip

where the characters $\chi_i: L^* \ra K^*$, $i=1,2$, are locally $L$-analytic.
We lift $\chi$ to a character of $\bP = \frP(L)$, and denote by $K_\chi$
the corresponding representation on $K$. Consider the locally
$L$-analytic representation $\Ind^\bG_\bP(K_\chi)$. For $t \in L^*$ close to $1$,
one has $\chi_i(t) = t^{c_i}$ with
$c_i = \frac{d}{dt}\chi_i(1)\in K$, $i=1,2$. The representation of $\frp$
on the dual space $(K_\chi)' \simeq K$ is given by the linear form

$$\left(\begin{array}{cc}
t_1 & \ast \\
0   & t_2
\end{array} \right) \mapsto - (c_1 t_1 + c_2 t_2) \,.$$

\medskip

By Thm. \ref{BGG} and Thm. \ref{main thm} we conclude that $\Ind^\bG_\bP(K_\chi)$ is
topologically irreducible if $-(c_1-c_2) \notin \bbZ_{\ge 0}$.
This was shown also in \cite{KS}, cf. Thm. 3.1.6.
\end{example}

\bigskip

\subsection{Restriction of Scalars} Let $\tilde{L}$
be a finite extension of $L$ and denote by $\Gamma$ the set of $L$-embeddings  of $\tilde{L}$ into $K.$
Suppose that the cardinality of $\Gamma$  is equal
to $[\tilde{L}:L]$. Let $\tilde{\frG}$ be a connected split
reductive group over $\tilde{L}$, $\tilde{\frS}$ a maximal split
torus and $\tilde{\frP}$ a Borel subgroup containing
$\tilde{\frS}$. Then we consider the reductive group $\frG =
\Res^{\tilde L}_L \tilde{\frG}$ over $L$. Put $\frP = \Res^{\tilde
L}_L \tilde{\frP}$, $\frT = \Res^{\tilde L}_L \tilde{\frS}$, and

$$\bG = \frG(L) = \tilde{\frG}(\tilde{L}) \,, \,\,
\bP = \frP(L) = \tilde{\frP}(\tilde{L}) \,, \,\,
\bT = \frT(L) = \tilde{\frS}(\tilde{L}) \,.$$

\medskip

Consider a locally $L$-analytic character $\chi: \bP \ra K^*$
which factors through $\bT$, and let $d\chi \in \Hom_L(\frt,K)$ be
the derivative of $\chi$. Because $\frt = \frs
\otimes_L \tilde{L}$ and by our assumption on $K$ we have a
canonical $K$-linear isomorphism

$$\Hom_L(\frt,K) = \Hom_L(\frs,K) \otimes_K \Hom_L(\tilde{L},K) \,.$$

\medskip

The set $\Gamma$ is a basis for the set of $L$-vector space
homomorphism $\Hom_L(\tilde{L},K)$, and so we can write

$$d\chi = \sum_{\sigma \in \Gamma} d\chi_\sigma \otimes \sigma \,,$$

\medskip

with uniquely determined $L$-linear maps $d\chi_\sigma: \frs \ra K$.

\medskip

The proof of the following assertion is an easy exercise:

\begin{prop} (i) The map

$$\tilde{L} \otimes_L K \lra \prod_{\sigma \in \Gamma}K \,, \,\,
\xi \otimes \zeta \mapsto (\sigma(\xi)\zeta)_\sigma \,,$$

\medskip

induces an isomorphism $U(\frg) \otimes_L K \simeq \prod_{\sigma
\in \Gamma} U(\tilde{\frg}_K)$ such that the $U(\frg) \otimes_L
K$-module $U(\frg_K) \otimes_{U(\frp_K)} K_{-d\chi}$ becomes
isomorphic to the direct sum of the modules

$$U(\tilde{\frg}_K) \otimes_{U(\tilde{\frp}_K)} K_{-d\chi_\sigma} \,.$$

\medskip

(ii) The $U(\frg) \otimes_L K$-module $U(\frg_K)
\otimes_{U(\frp_K)} K_{-d\chi}$ is irreducible if and only if for
every $\sigma \in \Gamma$ the $U(\tilde{\frg})$-module
$U(\tilde{\frg}) \otimes_{U(\tilde{\frp})} K_{-d\chi_\sigma}$ is
irreducible. If the latter condition is fulfilled, then
$\Ind^\bG_\bP(K_\chi)$ is a topologically irreducible
representation.
\end{prop}

\medskip

\begin{example} With the notation introduced above we let
$\tilde{\frg} = \Lie(\GL_{2,\tilde{L}})$ and hence we consider
$\bG = \GL_2(\tilde{L})$ as a locally $L$-analytic group. Let
$\chi$ be a locally $L$-analytic character of

$$\bT = \left\{ \left( \begin{array}{cc}
t_1 & 0 \\
0   & t_2
\end{array} \right) \midc t_1,t_2 \in \tilde{L}^* \right\} \,.$$

\medskip

For $t_1,t_2 \in \tilde{L}^*$ sufficiently close to $1$ we can write

$$\chi \left(\begin{array}{cc}
t_1 & 0 \\
0   & t_2
\end{array} \right) =
\prod_{\sigma \in \Gamma}
\sigma(t_1)^{c_{1,\sigma}}\sigma(t_2)^{c_{2,\sigma}} \,.$$

\medskip

Then the representation $\Ind^\bG_\bP(K_\chi)$ is topologically
irreducible if for all $\sigma \in \Gamma$:

$$-(c_{1,\sigma} - c_{2,\sigma}) \notin \bbZ_{\ge 0} \,.$$

\medskip
\end{example}

\section{Appendix: completed distributions rings which are integral domains}\label{appendix}

In this section we consider $L$-uniform groups $H$ as defined in Remark \ref{L_uniform}, and we show that there is a sequence $r_m \in (\frac{1}{p},1) \cap p^\bbQ$, which tends to 1, such that the completed distribution rings $D_{r_m}(H,K)$ are integral domains.

\bigskip

\subsection{The case $H = \fro_L$} We begin by considering the case when $H = \fro_L$. We denote by $q$ the cardinality of the residue field of $L$, and by $e$ the ramification index of $L$ over $\Qp$. Fix a basis $\Zp$-basis $v_1=1, \ldots, v_n$ of $\fro_L$. Denote by $h_i = \delta_{v_i}$ the corresponding delta distribution in $D(\fro_L,K)$. For given $m \in \bbZ_{\ge 0}$ we denote by $\cF_m(\fro_L,\Cp)$ the space of all maps $\fro_L \ra \Cp$ which are $L$-rigid-analytic on the cosets of $p^m\fro_L$.

\bigskip

\begin{prop}\label{injectivity1} (i) Suppose $s \in (\frac{1}{p},1) \cap p^{\Qp}$ has the property that $p^{-1/(p-1)-1/e} \le s^{\vep_p} < p^{-1/(p-1)}$. Then there is a unique continuous $K$-linear map

$$\phi_s: D_s(\fro_L,K) \lra  \cF_0(\fro_L,\Cp)'_b = \Hom^{cont}_K(\cF_0(\fro_L,\Cp),\Cp)_b$$

\medskip

which makes the diagram

$$\begin{array}{ccc}
D_s(\fro_L,K) & \stackrel{\phi_s}{\lra} & \cF_0(\fro_L,\Cp)'_b\\
\lhookuparrow & & \lhookuparrow\\
D(\fro_L,K) & \lra & C^{an}(\fro_L,\Cp)'_b
\end{array}$$

\medskip

commutative. The vertical maps being the canonical ones. Moreover, $\phi_s$ is injective.

\medskip

(ii) Suppose $r \in (\frac{1}{p},1) \cap p^{\bbQ}$ is such that, for some $m \in \bbN_0$, the number $s = r^{p^m}$ has the property that $s > \frac{1}{p}$ and  $p^{-1/(p-1)-1/eq^{em}} \le s^{\vep_p} < p^{-1/(p-1)}$. Then there is a unique continuous $K$-linear map

$$\phi_r: D_r(\fro_L,K) \lra  \cF_m(\fro_L,\Cp)'_b = \Hom^{cont}_K(\cF_m(\fro_L,\Cp),\Cp)_b$$

\medskip

which makes the diagram

$$\begin{array}{ccc}
D_s(\fro_L,K) & \stackrel{\phi_s}{\lra} & \cF_0(\fro_L,\Cp)'_b\\
\uparrow & & \lhookuparrow\\
D_r(\fro_L,K) & \stackrel{\phi_r}{\lra} & \cF_m(\fro_L,\Cp)'_b
\end{array}$$

\medskip

commutative. The vertical maps being the canonical ones. Moreover, $\phi_r$ is injective. In particular, the canonical map $D_r(\fro_L,K) \ra D_s(\fro_L,K)$ on the left of the diagram is injective.

\medskip

\end{prop}

{\it Proof.} {\it Step 1.} Before proving (i) and (ii) we will carry out some preliminary considerations related to functions in $\cF_m(\fro_L,\Cp)$. To do this we will heavily rely on the paper \cite{ST4}, especially sections 3 and 4. By \cite{ST4}, Prop. 4.5, any function $f \in \cF_m(\fro_L,\Cp)$ has a generalized Mahler expansion of the form

$$f(z) = \sum_{\beta \ge 0} c_\beta P_\beta(z\Omega)$$

\medskip

with $|c_\beta| p^{\beta/eq^{em-1}(q-1)} \ra 0$. The polynomials $P_\beta(Y)$ are those defined in \cite{ST4}, Def. 4.1. $\Omega = \Omega_{t_0'}$ is the period of the Lubin-Tate group $\cG$ in \cite{ST4}. Consider a distribution of the form

$${\bf b}_{(m)}^\gamma := (h_1^{p^m}-1)^\gamma \,.$$

\medskip

By \cite{ST4}, Lemma 4.6 (2), we have

$${\bf b}_{(m)}^\gamma\Big(P_\beta(z\Omega)\Big) = \Big\{(F_{p^mt'_0})^\gamma,P_\beta(- \cdot \Omega)\Big\}$$

\medskip

and by (9) of that Lemma

$$\Big\{(F_{p^mt'_0})^k,P_\ell(- \cdot \Omega)\Big\}
 = \frac{1}{\ell!} \frac{d^\ell(F_{p^mt'_0})^k}{dZ^\ell} (0) \,,$$

\medskip

where the power series $F_{at'_0}(Z) \in \fro_{\Cp}[[Z]]$ is defined in section 3 of \cite{ST4}. It follows from the formula before Prop. 3.1 in \cite{ST4} that

$$F_{p^mt'_0}(Z) = F_{t'_0}([p^m]_\cG(Z)) \,.$$

\medskip

Put $\rho_m = p^{-1/eq^{em}(q-1)}$. By \cite{ST4}, Lemma 3.2, the map $[p^m]_\cG$ maps
$\bB(\rho_m)$ to $\bB(\rho_0)$. Here, $\bB(\rho)$ denotes the affinoid disk of radius $\rho$. As $\rho_0 = p^{-1/e(q-1)}$ we deduce from \cite{ST4}, Lemma 3.4 (c), and the proof given there, that $F_{t'_0}$ maps $\bB(\rho_0)$ into $\bB(\rho_0|\Omega_{t'_0}|) = \bB(p^{-1/(p-1)})$. This implies that

\begin{numequation}\label{shrinking}
|Z| \le \rho_m \hskip5pt \Rightarrow \hskip5pt |F_{p^mt'_0}(Z)| \le p^{-1/(p-1)} \,.
\end{numequation}

\medskip

Next we write

$$F_{p^mt'_0}(Z) = \lambda_1Z + \lambda_2Z^2 + \lambda_3Z^3 + \ldots \in \fro_{\Cp}[[Z]] \,.$$

\medskip

It follows from (\ref{shrinking}) that $|\lambda_\ell|\rho_m^\ell \le p^{-1/(p-1)}$ for all $\ell \ge 1$, i.e.,

\begin{numequation}\label{lambda}
|\lambda_\ell| \le p^{-\frac{1}{p-1} + \frac{\ell}{eq^{em}(q-1)}} \,.
\end{numequation}

\medskip

Write

$$F_{p^mt'_0}(Z)^k =  \sum_{\ell \ge k} \Big(\sum_{
\begin{array}{c}\ell_1 + \ldots + \ell_k = \ell\\
\ell_i \ge 1
\end{array}} \lambda_{\ell_1} \cdot \ldots \cdot \lambda_{\ell_k}\Big) Z^\ell = \sum_{\ell \ge 0} \lambda_{k,\ell}Z^\ell \,.$$

\medskip

Note that $\lambda_{k,\ell}$ vanishes for $k > \ell$. Using (\ref{lambda}) this shows that

\begin{numequation}\label{estim}
\Big|\Big\{(F_{p^mt'_0})^k,P_\ell(- \cdot \Omega)\Big\}\Big|
 = |\lambda_{k,\ell}| \le p^{-\frac{k}{p-1} + \frac{\ell}{eq^{em}(q-1)}} \hskip8pt .
\end{numequation}

\medskip

Next we consider a (formal) sum of distributions

$$\mu = \sum_{\gamma \ge 0} d_\gamma {\bf b}_{(m)}^\gamma$$

\medskip

with $|d_\gamma|s^{\vep_p \gamma} \ra 0$ as $\gamma \ra \infty$. Applying $\mu$ to a function

$$f = \sum_{\beta \ge 0} c_\beta P_\beta(z\Omega) \in \cF_m(\fro_L,\Cp)$$

\medskip

gives the formal sum

\begin{numequation}\label{formal}
\mu(f) = \sum_{\beta \ge 0, \gamma \ge 0} d_\gamma c_\beta \Big\{(F_{p^mt'_0})^\gamma,P_\beta(- \cdot \Omega)\Big\}
\end{numequation}

\medskip

As we have seen above, the term

\begin{numequation}\label{term}
d_\gamma c_\beta \Big\{(F_{p^mt'_0})^\gamma,P_\beta(- \cdot \Omega)\Big\}
\end{numequation}

\medskip

vanishes if $\gamma > \beta$. In particular, we only need to consider the case when $\gamma \le \beta$. It follows from (\ref{estim}) that the absolute value of (\ref{term}) can be bounded from above as follows

\begin{numequation}\label{est2}
\begin{array}{rl}
& \left|d_\gamma c_\beta \Big\{(F_{p^mt'_0})^\gamma,P_\beta(- \cdot \Omega)\Big\}\right| \\
&\\
\le & |d_\gamma||c_\beta| p^{-\frac{\gamma}{p-1} + \frac{\beta}{eq^{em}(q-1)}}\\
&\\
= & \big(|d_\gamma|s^{\vep_p \gamma}\big) \cdot \big(|c_\beta| p^{\frac{\beta q}{eq^{em}(q-1)}}\big) \cdot \left(\frac{p^{-1/(p-1)}}{s^{\vep_p}}\right)^{\gamma} \cdot p^{-\frac{\beta}{eq^{em}}}\\
&\\
 \le & \big(|d_\gamma|s^{\vep_p \gamma}\big) \cdot \big(|c_\beta| p^{\frac{\beta q}{eq^{em}(q-1)}}\big) \cdot \left(\frac{p^{-1/(p-1)}}{s^{\vep_p}}\right)^{\beta} \cdot p^{-\frac{\beta}{eq^{em}}}\\
&\\
 = & \big(|d_\gamma|s^{\vep_p \gamma}\big) \cdot \big(|c_\beta| p^{\frac{\beta q}{eq^{em}(q-1)}}\big) \cdot \left(\frac{p^{-1/(p-1) - 1/eq^{em}}}{s^{\vep_p}} \right)^{\beta}
\end{array}
\end{numequation}

\medskip

By our assumptions on $d_\gamma$, $c_\beta$, and because $p^{-1/(p-1) - 1/eq^{em}} \le s^{\vep_p}$, we see that the sum (\ref{formal}) converges.

\medskip

{\it Step 2: proof of (i).} We apply \cite{Sch}, Prop. 5.9, to the ring $D_s(\fro_L,K)$. Hence we know that any $\mu \in D_s(\fro_L,K)$ has a unique expansion

$$\mu = \sum_{\gamma \ge 0} d_\gamma {\bf b}^\gamma \,,$$

\medskip

where

$${\bf b}^\gamma = (h_1-1)^\gamma \,,$$

\medskip

with $|d_\gamma|s^{\vep_p \gamma} \ra 0$, and the norm $\|\mu\|_s$ can be calculated as

$$\|\mu\|_s = \sup_\gamma |d_\gamma| s^{\vep_p \gamma} \,.$$

\medskip

From step 1 (with $m=0$) we know that any such $\mu$ defines a linear form $\phi_s(\mu)$ on $\cF_0(\fro_L,\Cp)$ which is obviously continuous. Because

\begin{numequation}\label{Kronecker_relation}
{\bf b}^\gamma \left( {z \choose \beta} \right) = \left\{\begin{array}{rcl}
1&,& \gamma = \delta\\
0 & ,& \gamma \neq \delta
\end{array} \right. \,,
\end{numequation}

\medskip

the map $D_s(\fro_L,K) \ra \cF_0(\fro_L,\Cp)'$ is injective. We deduce from (\ref{formal}) and (\ref{est2}) that $|\phi_s(\mu)(f)| \le \|\mu\|_s \|f\|$. This proves that $\phi_s$ is continuous. Further $\phi_s$ is uniquely determined by the commutativity of the diagram, because the image of $D(\fro_L,K)$ in $D_s(\fro_L,K)$ is dense. That the diagram is commutative follows from the fact that the group algebra $K[\fro_L]$ is dense in $D(\fro_L,K)$, and because any element $h_i -1$ has a convergent expansion in $D_s(\fro_L,K)$:

$$h_i -1 \hskip6pt = \hskip6pt ''((h_1-1)+1)^{v_i} -1'' \hskip6pt = \hskip6pt \sum_{\nu = 1}^\infty {v_i \choose \nu}(h_1-1)^\nu \,.$$

\medskip

For a function $f$ which is a polynomial of the form ${z \choose \beta}$, it follows from (\ref{Kronecker_relation}) that

$$\phi_s\left(h_i-1\right)\left(f\right) = \phi_s\left(\sum_{\nu = 1}^\infty {v_i \choose \nu}(h_1-1)^\nu\right)\left(f\right) = f(v_i) - f(0)$$

\medskip

Approximating any $L$-rigid-analytic function $f$ on $\fro_L$ by polynomials, this holds true for all elements of $\cF_0(\fro_L,\Cp)$. This proves the first assertion.

\medskip

{\it Step 3: proof of (ii).} The subgroup $p^m \fro_L \sub \fro_L$ is uniform pro-$p$, and we can thus consider the canonical $p$-valuation on this subgroup, and the associated norm $\|\cdot\|_s$ which we denote for clarity by $\|\cdot\|_s^{(m)}$. By \cite{Sch}, Lemma 7.4, the norm $\|\cdot\|_r$ on $D(\fro_L,K)$, when restricted to $D(p^m \fro_L,K) \sub D(\fro_L,K)$ induces the same topology on $D(p^m \fro_L,K)$ as the norm $\|\cdot\|_s^{(m)}$. Furthermore, $D_r(\fro_L,K)$ is a finite and free (left or right) module over the subring $D_s(p^m \fro_L,K)$, on a basis any set $\cR$ of coset representatives for $\fro_L/p^m \fro_L$:

$$D_r(\fro_L,K) = \bigoplus_{a \in \cR} (\delta_a)_r D_s(p^m \fro_L,K) \,.$$

\medskip

Here, $(\delta_a)_r$ denotes the image of the delta distribution $\delta_a \in D(\fro_L,K)$ in $D_r(\fro_L,K)$. Because of our assumptions on $s$ we can apply \cite{Sch}, Prop. 5.9, which tells us that any element $\mu \in D_s(p^m \fro_L,K)$ has a unique expansion

$$\mu = \sum_{\gamma \ge 0} d_\gamma {\bf b}_{(m)}^\gamma \,,$$

\medskip

where

$${\bf b}_{(m)}^\gamma = (h_1^{p^m}-1)^\gamma \,,$$

\medskip

and the norm $\|\mu\|_s^{(m)}$ can be calculated as

$$\|\mu\|_s^{(m)} = \sup_\gamma |d_\gamma|s^{\vep_p \gamma} \,.$$

\medskip

As we have seen in step 1, any $\mu \in D_s(p^m \fro_L,K)$ induces a continuous linear form $\phi_r(\mu)$ on $\cF_m(\fro_L,\Cp)$. By step 2 we have $\phi(\delta_a)(f) = f(a)$ for any $a \in p^m\fro_L$ and $f \in \cF_m(\fro_L,\Cp)$. In particular, $\phi(\mu)$ vanishes on all functions $f$ which vanish identically on $p^m\fro_L$, because the group ring $K[p^m\fro_L]$ is dense in $D_s(p^m\fro_L,K)$. For a product $(\delta_a)_r \mu$ with $a \in \cR$ and $\mu \in D_s(p^m \fro_L,K) \sub D_r(\fro_L,K)$, and $f \in \cF_m(\fro_L,\Cp)$ define

$$\left((\delta_a)_r \mu\right)(f) = \mu \left[\fro_L \ra \Cp \,, \, z \mapsto f(a + z)\right]$$

\medskip

Of course, for any $a \in \fro_L$ and any $f \in \cF_m(\fro_L,\Cp)$, the function $z \mapsto f(a + z)$ is again in $\cF_m(\fro_L,\Cp)$. Thus we have a well-defined map

$$\phi_r: D_r(\fro_L,K) = \bigoplus_{a \in \cR} (\delta_a)_r D_s(p^m \fro_L,K) \lra \cF_m(\fro_L,\Cp)'_b$$

\medskip

That this map is continuous is easily seen. To see that it is injective, fix $a_0 \in \cR$, and consider the subspace of those $f \in \cF_m(\fro_L,\Cp)$ which vanish outside $a_0 + p^m \fro_L$. Then, if an element

$$\lambda = \sum_{a \in \cR} (\delta_a)_r \mu_a \in D_r(\fro_L,K) \,,$$

\medskip

with all $\mu_a \in D_s(p^m \fro_L,K)$, is such that $\phi_r(\lambda) = 0$, then, for $f$ in the subspace just mentioned,

$$0 = \sum_{a \in \cR} \mu_a \left[z \mapsto f(a + z)\right] = \mu_{a_0}\left[z \mapsto f(a_0 + z)\right]$$

\medskip

But this means that $\mu_{a_0}$ vanishes on all $L$-rigid-analytic functions supported on $p^m \fro_L$. Because $\phi_s$ for $p^m \fro_L$ is injective (by step 2), this means that $\mu_{a_0} = 0$. Therefore, $\phi_r$ is injective as well. \qed

\bigskip

\subsection{Reduction to the case $H = \fro_L$.}

\begin{prop}\label{injectivity2} Let $H$ be a locally $L$-analytic group which is $L$-uniform (cf. Remark \ref{L_uniform}). We define the norms $\bar{q}_r$ using the canonical $p$-valuation on $H$ (cf. \ref{norms}).

\medskip

(i) Let $r \in (\frac{1}{p}) \cap p^\bbQ$ and $s = r^{p^m}$ be as in \ref{injectivity1} (ii). Then the canonical map $D_r(H,K) \ra D_s(H,K)$ is injective.

\medskip

(ii) There is a sequence of numbers $r_m \in (\frac{1}{p},1) \cap p^\bbQ$, which tends to 1, such that the completed distribution rings $D_{r_m}(H,K)$ are integral domains.
\end{prop}

{\it Proof.} (i) Let $\Lambda$ be as in Lemma \ref{L_uniformgroups}. By \cite{DDMS}, 9.10, we see that $\Lambda$ is a uniform $\Zp$-Lie algebra: $[\Lambda,\Lambda] \sub p^{\vep_p}\Lambda$ . By \cite{DDMS}, 9.8, the Baker-Campbell-Hausdorff series $\Phi(X,Y) = \exp_H^{-1}(\exp_H(X)\exp_H(Y))$ converges on $\Lambda \times \Lambda$ and maps $\Lambda \times \Lambda$ into $\Lambda$. In fact, we can identify $H = \exp_H(\Lambda)$ with $\Lambda$, the group multiplication being given by $\Phi(\cdot, \cdot)$. The members of the lower $p$-series are then $\exp_H(p^k\Lambda)$, $k \ge 0$. Let $(v_1, \ldots, v_n)$ be a $\Zp$-basis of $\fro_L$, and let $(\frx_j)_{j=1}^d$ be an $\fro_L$-basis of $\Lambda$. Then  $(v_i \frx_j)_{1 \le i \le n, 1 \le j \le d}$ is a $\Zp$-basis for $\Lambda$. The map $\psi: \fro_L^d \stackrel{\sim}{\lra} H$ defined by sending

$$(z_j)_{j=1}^d = (a_{1,j}v_1 + \ldots + a_{n,j}v_n)_{j=1}^d \in \fro_L^d$$

\medskip

to

$$\begin{array}{l}
\exp_H(v_1 \frx_1)^{a_{1,1}} \cdot \ldots \cdot  \exp_H(v_n \frx_1)^{a_{n,1}} \cdot \ldots \cdot \exp_H(v_1 \frx_d)^{a_{1,d}} \cdot \ldots \cdot\exp_H(v_n \frx_d)^{a_{n,d}}  \\
\\
=  \exp_H(z_1\frx_1) \cdot \ldots \cdot \exp_H(z_d\frx_d)\\
\end{array}$$

\medskip

is then a homeomorphism, cf. Lemma \ref{L_uniformgroups} (3). It is shown in \cite{ST3}, Prop. A.3, that

$$D(\fro_L^d,K) =  D(\fro_L,K) \widehat{\otimes}_{K,\pi} \ldots \widehat{\otimes}_{K,\pi} D(\fro_L,K)$$

\medskip

(With $d$ factors on the right.) Here, the completed topological tensor product is the completion of the ordinary tensor product with respect to the projective tensor product topology (indicated by $\pi$). In \cite{ST3}, Prop. A.3, it is the inductive tensor product topology that is used. However, by \cite{S}, 17.6, the inductive and the projective tensor product topologies coincide in this case, as the distribution algebras are Fr\'echet spaces. Using the homeomorphism $\psi$ we get an isomorphism of topological vector spaces

$$\psi^\ast: C^{an}_L(H,K) \stackrel{\simeq}{\lra} C^{an}_L(\fro_L^d,K)$$

\medskip

which gives, after dualizing, an isomorphism of topological vector spaces

$$(\psi^\ast)':  D(\fro_L,K) \widehat{\otimes}_{K,\pi} \ldots \widehat{\otimes}_{K,\pi} D(\fro_L,K) =  D(\fro_L^d,K) \stackrel{\simeq}{\lra} D(H,K)$$

\medskip

Tracing through the definition of the norms $\bar{q}_r$ on $D(H,K)$ as constructed in \ref{norms} shows that $\bar{q}_r$ corresponds to the norm $\bar{q}_r^{\fro_L} \otimes \ldots \otimes \bar{q}_r^{\fro_L}$ on

$$D(\fro_L,K) \widehat{\otimes}_{K,\pi} \ldots \widehat{\otimes}_{K,\pi} D(\fro_L,K)$$

\medskip

induced by the norms $\bar{q}_r^{\fro_L}$ on each factor $D(\fro_L,K)$. Here again, the norms $\bar{q}_r^{\fro_L}$ are defined as in \ref{norms}, using the canonical $p$-valuation on the uniform group $\fro_L$. The completion of

$$D(\fro_L,K) \widehat{\otimes}_{K,\pi} \ldots \widehat{\otimes}_{K,\pi} D(\fro_L,K)$$

\medskip

with respect to the induced norm $\bar{q}_r^{\fro_L} \otimes \ldots \otimes \bar{q}_r^{\fro_L}$ is then

$$D_r(\fro_L,K) \widehat{\otimes}_{K,\pi} \ldots \widehat{\otimes}_{K,\pi} D_r(\fro_L,K)$$

\medskip

It follows from \cite{Em}, 1.1.27, that if the map $D_r(\fro_L,K) \ra D_s(\fro_L,K)$ is injective, then so is the (iterated) tensor product of this map with itself. By Prop. \ref{injectivity1} we know that this map is injective if $r$ and $s$ are in relation as stated in Prop. \ref{injectivity1}. Therefore we have proved that the canonical map $D_r(H,K) \ra D_s(H,K)$ is injective (for $r$ and $s$ as in Prop. \ref{injectivity1}).

\medskip

(ii) It follows from \cite{Sch}, Prop. 5.6, that if $s> \frac{1}{p}$ and $s^{\vep_p} < p^{-1/(p-1)}$, then $D_s(H,K)$ is an integral domain. Hence $D_r(H,K)$ is an integral domain. It is easily seen, that one can find an increasing sequence of numbers $r_m \in (\frac{1}{p},1) \cap p^\bbQ$ which tends to 1, and such that $s_m = r_m^{p^m}$ has the property that $s_m> \frac{1}{p}$ and $p^{-1/(p-1)-1/eq^{em}} \le s_m^{\vep_p} < p^{-1/(p-1)}$. \qed

\end{document}